\documentclass[10pt,a4paper]{article}
\usepackage[french]{babel}
\usepackage[latin1]{inputenc}

\usepackage[centertags]{amsmath}
\usepackage{amsfonts}
\usepackage{amssymb}
\usepackage{amsthm}
\usepackage{newlfont}

\usepackage{fancyhdr}
 \fancypagestyle{plain}{\fancyhead{}
}

\newtheorem{prop}{Proposition}
\newtheorem{defi}{Définition}
\newenvironment{preuve}{\textsl{Preuve.}}{$\blacksquare\\$}
\newenvironment{rem}{\textbf{Remarque.}}{\\}
\newenvironment{ex}{\textbf{Exemple$^{}$}}{}

\newtheorem{lem}{Lemme}
\newtheorem{coro}{Corollaire}

\newcommand{\derivnieme}[3]{\frac{d^{#3}#1}{d#2^{#3}}}
\newcommand{\derivpartnieme}[3]{\frac{\partial^{#3}#1}{\partial#2^{#3}}}

\newcommand{\complex}{\mathbb C}
\newcommand{\integer}{\mathbb N}
\newcommand{\rinteger}{\mathbb Z}
\newcommand{\real}{\mathbb R}
\newcommand{\seriesfactform}{\widehat{\mathcal{O}}_{fact}}
\newcommand{\corpsseriesfactform}{\widehat{\mathcal{M}}_{fact}}
\newcommand{\seriesfact}{\mathcal{O}_{fact}}
\newcommand{\corpsseriesfact}{\mathcal{M}_{fact}}
\newtheorem{theo}{Théorème}
\newcommand{\seriesretrofact}{\mathcal{O}_{r\acute{e}tro-fact}}
\newcommand{\seriesretrofactform}{\widehat{\mathcal{O}}_{r\acute{e}tro-fact}}
\newcommand{\corpsseriesretrofact}{\mathcal{M}_{r\acute{e}tro-fact}}
\newcommand{\corpsseriesretrofactform}{\widehat{\mathcal{M}}_{r\acute{e}tro-fact}}

\newcommand{\io}{]0,h_0[}
\begin{document}


\title{\Huge Classification rationnelle et confluence des systèmes aux différences singuliers réguliers
}
\author{par Julien ROQUES \\
$_{}$\\
\small \texttt{E-mail : roques@picard.ups-tlse.fr}\\
\small Laboratoire Emile Picard, Université Paul Sabatier\\
\small 118 route de Narbonne, 31062 Toulouse Cedex, France.\\
\small \textbf{Mots-clés} \textbf{:} Equations
\small aux différences, matrice de connexion,\\
\small équations différentielles, monodromie. \\
\small \textbf{Classification
 \small math.} \textbf{:} 39A10, 34M35, 34M40, 65Q05.\\}

\date{}
\maketitle{}


\scriptsize \tableofcontents \small
\newpage
\hrule
  \abstract{En choisissant des ``caractères'' et des ``logarithmes'', méromorphes sur $\complex,$ construits à l'aide de la fonction
  Gam\-ma d'Euler, et en utilisant des séries de factorielles convergentes, nous sommes en mesure,
  dans une première partie, de donner une ``forme normale'' pour les solutions
  d'un système aux différences singulier
  régulier. Nous pouvons alors définir une matrice de connexion d'un tel système.
  Nous étudions ensuite, suivant une idée de G.D. Birkhoff, le lien de celles-ci avec le problème de la classification rationnelle des systèmes. Dans une deuxième partie, nous nous intéressons à la confluence
  des systèmes aux différences fuchsiens vers les systèmes différentiels. Nous montrons en particulier comment, sous certaines hypothèses naturelles, on peut reconstituer les monodromies locales d'un système différentiel limite à partir des matrices méromorphes de connexion des déformations considérées. Le point central,
  qui distingue
  en profondeur les systèmes aux différences singuliers réguliers de leurs homonymes différentiels ou aux $q$-différences et qui rend
  leur étude plus complexe, est la nécessaire utilisation de séries
  de factorielles (qui peuvent diverger en tant que séries de puissances). Une version de cet article est à paraître aux Annales de l'Institut Fourier.
   \begin{center}
 \textbf{Abstract} \\
   \end{center}

 By using meromorphic ``characters'' and ``logarithms'' built up from Euler's Gamma function, and
by using convergent factorial series, we will give, in a first pat, a ``normal form'' to the solutions of a singular
regular system. It will enable us to define a connexion
matrix for a regular singular system. Following one of Birkhoff's
idea, we will then study its link with the problem of rational
classification of systems. In a second part, we will be interested in the
confluence of fuchsian difference systems to differential systems.
We will show more particularly how we can get, under some natural
hypotheses, the local monodromies of a limit differential system from the
connection matrices of the deformation that we consider.
The use of factorial series (which can diverge as power series) distinguish regular
singular difference systems from their differential and
$q$-difference analogues and make their study more difficult.}\\
\hrule

\normalsize
\section{Introduction.}

Un système aux différences :
\begin{equation} \label{le systeme}
  Y(x-1)=A(x)Y(x),\ A \in
Gl_{n}(\complex(x))
\end{equation}
est dit fuchsien si $A$ est holomorphe à l'$\infty$ et si
$A(\infty)=I.$ Il est dit singulier
régulier s'il peut être transformé en un système fuchsien à l'aide
d'une transformation de jauge rationnelle. Rappelons que le
système déduit de (\ref{le systeme}) par la
transformation de jauge $R$ ($R \in Gl_{n}(K)$ où $K$ est un extension de corps de $\complex(x)$)
est celui obtenu en posant
$Y=RZ$ :
\begin{equation*}
  Z(x-1)=\left[R(x-1)^{-1}A(x)R(x)\right]Z(x).
\end{equation*}
Désignant par $\delta_{-1}$ l'opérateur dont l'action sur une
fonction $Y$ est définie par la formule :
\begin{equation*}
\delta_{-1}Y(x)=(x-1)(Y(x)-Y(x-1)),
\end{equation*}
un système fuchsien peut s'écrire sous la forme :
\begin{equation*}
  \delta_{-1}Y=AY
\end{equation*}
où $A$ est holomorphe à l'infini.

Un système fuchsien non résonnant (\textit{i.e.} tel que deux valeurs
propres de $A(\infty)$ ne diffèrent pas par un entier relatif non
nul) admet une unique solution formelle de la forme :
\begin{equation*}
  (I+\sum_{s=1}^{+\infty}\widehat{Y}_s x^{-s})x^{K}.
\end{equation*}
La série $I+\sum_{s=1}^{+\infty}\widehat{Y}_s x^{-s}$ n'est en
général pas convergente; c'est la différence fondamentale des
systèmes aux différences singuliers réguliers avec leurs analogues
aux $q$-différences et différentiels. Dans \cite{vds}, M. Van der
Put et M. Singer démontrent qu'il existe deux solutions fondamentales
asymptotiques à $\widehat{Y}$ dans un demi-plan gauche pour l'une
et droit pour l'autre.
Notons que Birkhoff avait déjà prouvé des résultats analogues dans un cas
``irrégulier'' (voir \cite{Birkhoff autre}).

D'autres auteurs ont remarqué l'intérêt, pour l'étude de ce
pro\-blème, des séries de factorielles (on peut citer N.-E.
Nörlund \cite{Norlund series d interp} et plus tard Fitzpatrick et
Grimm \cite{fitz grimm}, ou encore W.A. Jr. Harris \cite{harris}).
Rappelons qu'il s'agit de séries de la forme :

\begin{equation*}
  \sum_{s=0}^{+\infty}A_{s}x^{-[n]}
\end{equation*}
où :
\begin{equation*}
x^{-[n]}=\frac{1}{x(x+1)\cdots (x+n-1)}.
\end{equation*}
Rappelons également qu'en faisant le changement $x\leftarrow -x$
dans la définition des séries de factorielles, on obtient la
notion de série de rétro-factorielles. On peut démontrer que la série $I+\sum_{s=1}^{+\infty}\widehat{Y}_{s}x^{-s}$ est d\'{e}\-velop\-pable en série de
factorielles (resp. rétro-factorielles) convergente dans un
demi-plan droit (resp. gauche) et tangente à $I$ en $+\infty$
(resp. $-\infty$); voir par exemple \cite{harris}. Il convient également de remarquer que, d'un
point de vue combinatoire,
les $x^{-[n]}$ sont des vecteurs propres de l'opérateur $\delta_{-1}.$

Une dernière approche consiste à se ramener, par une
transformation de jauge développable en série de factorielles
convergente, tangente à $I$ en $+\infty,$ à l'équation
$Y(x-1)=(I-\frac{A(\infty)}{x-1})Y(x).$ Pour résoudre cette
équation, on se ramène à deux familles d'équations de bases :
celle des caractères et celle des logarithmes (elles
correspondent respectivement à la partie semi-simple et à la
partie nilpotente de $A(\infty)$), celles-ci peuvent être résolues
par des fonctions (uniformes) méromorphes sur $\complex$
construites à l'aide de la fonction Gamma d'Euler (il reste un arbitraire
dans le choix de telles solutions). Par exemple, A. Duval utilise
cette approche, dans \cite{duvconfqfuchsien}, pour étudier la
confluence des systèmes aux \textit{q}-différences fuchsiens vers
les systèmes aux différences fuchsiens. Nous adoptons aussi
cette démarche dans cet article. On peut opérer de manière
symétrique en $-\infty$ par le changement de variable $x\leftarrow
-x $.

Ainsi dans le cas non résonnant on a deux solutions fondamentales
``canoniques'' : l'une est attachée à $+\infty,$ l'autre à
$-\infty$. En revanche dans le cas résonnant ou plus généralement
dans le cas singulier régulier, il n'y a pas \textit{a priori} de
solution canonique, contrairement au cas différentiel. Se pose
donc le problème de donner une ``forme normale'' aux solutions.
Nous en proposons une, inspirée des travaux de J. Sauloy dans
\cite{sauloy systemes aux q diff singu reg }.

Cela nous permet de définir la matrice de connexion de Birkhoff d'un
système aux différences singulier régulier.
Celle-ci rend compte des relations linéaires entre les deux
systèmes fondamentaux de solutions évoqués ci-dessus. G.D.
Birkhoff a étudié son lien avec le problème de la classification
rationnelle des systèmes aux différences (voir \cite{Birkhoff
generalised RH pb} et \cite{Birkhoff autre}). Dans \cite{vds} M.
van der Put et M.F. Singer étudient la matrice de connexion de
manière approfondie dans le cas régulier et donnent quelques
énoncés dans le cas singulier régulier sans toujours détailler les
preuves. Nous réinterprétons ces énoncés dans notre approche et les prouvons en détail.\\

La deuxième partie de cet article est consacrée à la
confluence des systèmes aux différences fuchsiens vers les
systèmes différentiels : nous étudions la confluence des solutions
et celle des matrices de connexion. En particulier, nous
expliquons comment peuvent être calculées, à partir des matrices de connexion, les monodromies locales d'un système
différentiel limite.

L'étude de cette confluence est plus difficile
que celle, analogue, des systèmes aux $q$-différences vers les
systèmes différentiels effectuée par J. Sauloy dans \cite{sauloy
systemes aux q diff singu reg }, dont nous nous sommes inspirés, parce qu'on ne peut pas
employer des séries de $1/x$ convergentes, et qu'il faut les remplacer par
des séries convergentes de ($h$-)factorielles.

Notre approche repose sur une hypothèse naturelle relative à la
croissance des coefficients des séries de $h$-factorielles
qui définissent la déformation du système différentiel
considéré (hypothèses de type $(C,\lambda)$). Notons que ces
hypothèses sont vérifiées lorsqu'on est en présence de convergence
uniforme au voisinage de l'infini (voir le corollaire \ref{ca c
est un super resultat} en \ref{etude globale}), ce qui en fait des
hypothèses raisonnables, généralement vérifiées dans la pratique.

Signalons que dans \cite{krichever} I. Krichever présente une
approche différente de l'étude des systèmes aux différences
(techniquement, tout repose sur la résolution de
problèmes du type Riemann-Hilbert). Il propose en particulier,
pour certains systèmes aux différences, une notion de monodromies
locales. Il montre que dans certaines situations, ces monodromies
locales confluent vers les monodromies locales d'un système
différentiel limite. Les hypothèses sont toutefois assez restrictives.\\

Les problèmes que nous considérons (classification rationnelle et confluence) et la ``philosophie'' de leurs solutions sont
analogues à ceux de J. Sauloy dans
$\cite{sauloy systemes aux q diff singu reg }$. Ici cependant, l'apparition de
séries de $h$-factorielles rend les
situations rencontrées, et les moyens mis en oeuvre pour les
comprendre, plus complexes.

Notre article est également inspiré de travaux récents
de A. Duval; notamment des articles \cite{duvseriesfact} et \cite{duvconfqfuchsien}.\\

Mentionnons enfin que dans l'article en préparation \cite{duv roques} (voir également \cite{ma these}) nous envisagerons, en collaboration avec A. Duval, une généralisation des différents phénomènes de confluence.\\
\\
\textbf{Remerciements. } Ce travail fait partie d'une thèse
(\cite{ma these}) sous la direction de J.-P. Ramis; je lui adresse
toute ma reconnaissance pour son écoute et son aide si précieuses.
Je tiens également à exprimer toute ma gratitude envers J. Sauloy pour ne s'être jamais
économisé lors de nos nombreuses $q$-discussions. Enfin, je
remercie chaleureusement A. Duval pour son aide et pour l'intérêt
qu'elle a porté à ce travail.

 \section{Notations et terminologie.} \label{notations et terminologies}

On note $\tau_{-h}$ l'opérateur de translation de pas $-h:$
$\tau_{-h}\ y(x)=y(x-h)$.
On définit également les deux opérateurs $\Delta_{-h}$ et
$\delta_{-h}$ par :
\begin{equation*}
\Delta_{-h}\ y(x)=\frac{y(x)-\tau_{-h}\ y(x)}{h} \text{ et } \delta_{-h}\ y(x)=(x-h)\Delta_{-h}\ y(x).
\end{equation*}
Le système (aux différences) de pas $h$ défini par une fonction
matricielle $A$ est par définition le système suivant :
\begin{equation*}
\delta_{-h}Y=AY.
\end{equation*}

Les anneaux de séries de $h$-factorielles et de séries de
$h$-rétro-facto\-riel\-les sont respectivement
notés $\seriesfact^{(h)}$ et
$\seriesretrofact^{(h)};$ les anneaux de séries
formelles correspondants seront respectivement notés
$\seriesfactform^{(h)}$ et $\seriesretrofactform^{(h)}.$ Les corps
des fractions des anneaux précédents seront respectivement notés
$\corpsseriesfact^{(h)},$ $\corpsseriesretrofact^{(h)},$
$\corpsseriesfactform^{(h)}$ et $\corpsseriesretrofactform^{(h)}.$ Les sous-corps de $\corpsseriesfact^{(h)}$ et de  $\corpsseriesretrofact^{(h)}$ constitués de leurs éléments méromorphes sur $\complex$ tout entier seront respectivement notés $\corpsseriesfact^{(h)}(\complex)$ et $\corpsseriesretrofact^{(h)}(\complex).$ 
Nous définissons :
\begin{equation*}
x^{-[n]_{h}}=\frac{1}{x(x+h)\cdots (x+(n-1)h)},\ x^{[n]_{h}}=\frac{1}{ x^{-[n]_{h}}}
\end{equation*}
et :
\begin{equation*}
x^{-\{n\}_{h}}=\frac{1}{x(x-h)\cdots (x-(n-1)h)}, \ x^{\{n\}_{h}}=\frac{1}{x^{-\{n\}_{h}}};
\end{equation*}
ce sont des vecteurs propres de $\delta_{-h}.$
Lorsque nous n'envisagerons pas d'étudier des propriétés de
confluence, nous nous placerons dans le cas $h=1$ et nous
omettrons alors $h$ dans toutes nos notations (par exemple
$x^{-[n]}=x^{-[n]_{1}}$, \textit{etc}).

Nous traiterons, en vue de la confluence, de familles de systèmes aux
différences indexés par des pas $h>0.$ Nous introduisons quelques
définitions afin d'alléger les énoncés et de dégager les notions
importantes. On se donne $h_0>0.$\\

\begin{defi}
Soit $(C,\lambda)\in \real^{+2}.$ Donnons nous, pour tout $h\in \io,$ une fonction $A^{(h)}$ développable en série de $h$-factorielles : 
\begin{equation*}
  A^{(h)}(x)=\sum_{s=0}^{+\infty}A_{s}^{(h)}x^{-[s]_{h}} \in M_n(\seriesfact^{(h)}).
\end{equation*}
La famille $A^{(h)}$, $h \in \io$ est dite de type $(C,\lambda)$ s'il existe $h'_0 \in ]0,h_0[$ tel que pour tout $h\in ]0,h_0'[$  et pour tout $s \in \integer^*$ :
\begin{equation*}
  \|A^{(h)}_{s}\| \leq
C\lambda^{[s-1]_{h}}
\end{equation*}
et si la famille des $\|A^{(h)}_{0}\|$, $h \in \io$ est bornée.
\end{defi}

\begin{defi}
Une famille de systèmes de pas $h \in \io$ :
\begin{equation} \label{aaa}
\delta_{-h}Y=A^{(h)}Y
\end{equation}
sera dite de Fuchs ou fuchsienne $(C,\lambda)$ en $+\infty$ si :
\begin{itemize}
\item pour tout $h\in \io$, $A^{(h)} \in M_n(\seriesfact^{(h)})$,
\item la famille $A^{(h)}$, $h \in \io$ est de type $(C,\lambda).$
\end{itemize}
Elle sera dite algébrique de Fuchs ou fuchsienne $(C,\lambda)$ en
$+\infty$ si, pour tout $h \in \io,$ $A^{(h)} \in M_{n}(\complex(x))$ et
si cette famille est fuchsienne $(C,\lambda)$ en $+\infty$.

La famille de Fuchs $(C,\lambda)$ en $+\infty$ (\ref{aaa}) est dite non
résonnante si, pour tout $h\in \io$, deux éléments du spectre $Sp(A^{(h)}(+\infty))$ ne
diffèrent pas par un entier relatif non nul.
\end{defi}

Nous avons des notions analogues en $-\infty$ grâce au changement de variable $x \leftarrow -x.$

Enfin, le système obtenu après la transformation de jauge $F \in Gl_n(K)$ où $K$ est une extension de corps de $\complex(x)$, à partir du système $\delta_{-h}Y=AY$, est par définition le système $\delta_{-h}Y=A^FY$ où $A^F(x)=\frac{I-F(x-h)^{-1}(I-hA(x))F(x)}{h}$ (\textit{i.e.} celui obtenu après le changement de variable $Y \leftarrow FY$).

\section{Résolution.}

Dans cette première partie, on étudie des systèmes
algé\-briques (\textit{i.e.} à coefficients rationnels) de pas
$h=1$. On dira qu'un tel
système, $\delta_{-1}Y=AY$,  est fuchsien s'il est défini par un élément $A$ de
$M_{n}(\complex(x))$ holomorphe à l'infini. Il sera de plus dit
non résonnant si deux éléments de $Sp(A(\infty))$ ne diffèrent pas
par un entier relatif non nul. Un système (rationnel)
$\tau_{-1} Y=BY$ sera dit singulier régulier s'il peut être
transformé en un système fuchsien à l'aide d'une transformation de
jauge rationnelle.
\subsection{Systèmes aux différences fuchsiens non ré\-son\-nants.}

\subsubsection{Rappels.}

Rappelons les étapes de la résolution, en $+\infty,$ des systèmes
fuchsiens non résonnants;  pour plus de détails, nous renvoyons le
lecteur à \cite{duvseriesfact} et à \cite{duvconfqfuchsien}. Soit :
\begin{equation} \label{aa}
  \delta_{-1}Y=AY
\end{equation}
un système (algébrique) fuchsien non résonnant, on note
$A_{0}=A(\infty)$. On peut montrer qu'il existe une unique
transformation de jauge $F \in Gl_{n}(\seriesfact)$ tangente à
$I$ en $+\infty$
qui ramène au système, à matrice constante, défini par $A_0$ (\textit{i.e.} $A^F=A_0$).\\
\\
\begin{rem}
Notons que $F(x)$ et $F^{-1}(x)$ sont bornées pour $|x|$ assez grand dans tout demi-plan droit.
\end{rem}

Il suffit donc de savoir résoudre le système
constant défini par $A_0$. Donnons nous, pour tout $c\in \complex$, une famille $l^{(k)}_{c}, \ k\in \integer$
de fonctions uniformes et méromorphes sur $\complex$ qui satisfont :
\begin{equation*}
\delta_{-1}l^{(k)}_{c}=cl^{(k)}_{c}+l^{(k-1)}_{c}.
\end{equation*}
Posons $e_c^+=l^{(0)}_{c}$ et $l^{(k)}=l^{(k)}_{0}.$ Ces deux fonctions vérifient :
\begin{eqnarray*}
\delta_{-1}e_c^+&=&c e_c^+ \text{ (équations des caractères),}\\
\delta_{-1}l^{(k)}&=&l^{(k-1)} \text{ (équations des logarithmes).}
\end{eqnarray*}
Nous choisissons :
$$e_{c}^{+}(x)=\frac{\Gamma(x)}{\Gamma(x-c)},\ \
l^{(k)}_{c}(x)=\frac{1}{k!}\frac{\partial^{k}}{\partial
c'^{k}}_{|c'=c} e_{c'}^{+}(x).$$

Ecrivons une réduction de Jordan de $A_{0}$ :
\begin{equation*}
A_{0}=P diag(c_{1}I_{\mu_{1}}+N_{\mu_1},...,c_{m}I_{\mu_{m}}+N_{\mu_m})P^{-1}.
\end{equation*}
Alors, un système fondamental de solutions du système défini par
$A_{0}$ en $+\infty$ est donné par:

\begin{equation*}
e_{A_{0}}^{+}:=P diag(e_{c_{1}I_{\mu_{1}}+N_{\mu_1}}^{+},..,e_{c_{m}I_{\mu_{m}}+N_{\mu_m}}^{+})P^{-1}
\end{equation*}
avec :
\begin{equation*}
e_{cI_{\mu}+N_\mu}^{+}=\left( \begin{smallmatrix} l^{(0)}_{c}&l^{(1)}_{c}&...&l^{(\mu-1)}_{c}\\
                 \vdots&l^{(0)}_{c}&...&l^{(\mu-2)}_{c}\\
                 \vdots&&\ddots&\vdots\\
                 0&...&0&l^{(0)}_{c}
\end{smallmatrix} \right).
\end{equation*}
On montre sans peine que la solution obtenue est indépendante de la réduction de Jordan
choisie.

Un système fondamental de solutions (canonique) en $+\infty$ du système fuchsien non
réson\-nant (\ref{aa}) est donc donné par le produit d'une
transformation de jauge $F \in Gl_{n}(\seriesfact)$ tangente à
$I$ en $+\infty$ et de $e^{+}_{A(\infty)};$ on le note $e^{+}_{A}$.\\

Les mêmes constructions peuvent être réalisées en $-\infty$ par le
changement de variable $x\longleftarrow -x;$ on note $e^{-}_{A}$
la solution canonique ainsi obtenue. Elle s'écrit comme le produit
d'une fonction développable en série de rétro-factorielles
tangente à $I$ en $-\infty$ et de $e^{+}_{A(\infty)}(-x).$

\subsubsection{Comportement asymptotique des solutions canoniques
dans le cas non résonnant.}\label{sous section comportement
asymptotique}

Le résultat suivant est classique.
\begin{prop} \label{sol formelle}
Le système (\ref{aa}) admet une unique solution formelle de la
forme :
\begin{equation*}
  (I+\sum_{s=1}^{+\infty}\widehat{Y}_{s}x^{-s})x^{K}
\end{equation*}
avec $K \in M_{n}(\complex)$ et on a $K=A_{0}(=A(\infty)).$
\end{prop}

La fonction
$\frac{\Gamma(x)}{\Gamma(x-c)}x^{-c}$ est développable en série de factorielles et est tangente à $1$ en
$+\infty$ (voir l'appendice); nous en déduisons que, pour $k\geq 1,$ la fonction :
\begin{equation*}
   l^{(k)}_{c}(x)x^{-c}
  - l^{(k-1)}_{c}(x)  \log(x) x^{-c}
  +...+
  l^{(0)}_{c}(x) (-1)^{k}\frac{\log(x)^{k}}{k!} x^{-c}
\end{equation*}
est développable en série de factorielles et est tangente à $0$ en
$+\infty.$ Il en résulte que la fonction $e^{+}_{A_{0}}x^{-A_{0}}
$ est développable en série de factorielles et est tangente à $I$
en $+\infty.$

La solution canonique, en $+\infty,$ du système (\ref{aa}) est de
la forme $Fe^{+}_{A_{0}}$ avec $F \in Gl_n(\seriesfact)$ tangent à $I$
en $+\infty.$ On note $\widehat{F}$ la série formelle des $x^{-1}$
correspondant à $F$ et $\widehat{G}$ celle correspondant à
$e^{+}_{A_{0}}x^{-A_{0}}.$ Il est clair que
$\widehat{F}\widehat{G}x^{A_{0}}$ est une solution formelle du
système considéré : c'est celle exhibée à la proposition
\ref{sol formelle}. De plus, suivant un résultat de B. Malgrange dans \cite{seminaire malgrange}, une série de factorielles
convergente est asymptotique dans un demi-plan droit à son écriture en série
formelle\footnote{Rappelons que l'anneau des séries de
factorielles formelles est isomorphe à celui des séries en $1/x,$
par un isomorphisme canonique qui respecte la valuation.}.
Par conséquent
$e^{+}_{A}=Fe^{+}_{A_{0}}$ est asymptotique, dans un certain demi-plan droit $\{ z\in \complex \ | \ Re(z)>M\}$, à la solution formelle
de la proposition \ref{sol formelle}. Cela signifie que pour tout $N\in \integer,$ il existe $C_N \in \real^+$ tel que :
$$\|e^{+}_{A}x^{-A_{0}}-\sum_{s=0}^{N-1}\widehat{Y}_{s}x^{-s}\| \leq C_N |x|^{-N}$$
si $x \in \{ z\in \complex \ | \ Re(z)>M\}$ est de module suffisamment grand.
En utilisant l'équation fonctionnelle de $e^+_A$, on voit sans difficulté que ce développement est valable sur \textit{tout} demi-plan droit. Nous pouvons ainsi énoncer le

\begin{theo}
La solution $e^{+}_{A}$ est holomorphe et non dégénérée dans un
demi-plan droit, et asymptotique sur tout demi-plan droit à la
solution formelle de la proposition précédente.
\end{theo}

Naturellement, nous avons un énoncé analogue pour $e_{A}^{-}(x)$
dont nous laissons la formulation au lecteur.

\begin{theo}
Si une solution $\Psi$ du système (\ref{aa}) admet un
dévelop\-pement asymptotique dans un demi-plan droit de la forme
$\widehat{Z}x^{B_{0}}$ avec $\widehat{Z} \in
Gl_{n}(\complex((x^{-1})))$ et $B_{0} \in M_{n}(\complex)$ alors, il
existe $N \in Gl_{n}(\complex)$ tel que $\Psi = e^{+}_{A} N.$
\end{theo}

\begin{proof}
Le résultat étant indépendant de la norme choisie, on peut la supposer sous-multiplicative.
On note $\Phi=e^{+}_{A}$ et $\widehat{Y}x^{A_{0}}$ son
développement asymptotique sur un (sur tout) demi-plan droit. Si $n\in \integer,$ on désigne par $\widehat{Y}_{|n}$ la $n$-ème somme partielle de $\widehat{Y}$; notations similaires pour $\widehat{Z}.$ Soit
$\omega$ la fonction méromorphe sur $\complex$ et 1-périodique définie par $\omega=\Phi^{-1}\Psi.$ Donnons nous $N\in \integer$ tel que, pour tout $x\in \complex$ avec $Re(x)\geq 1$, on ait :
$$\max \{\|x^{A_0}\| \|x^{-B_0}\|,\|x^{-A_0}\| \|x^{B_0}\| \} \leq |x|^N .$$
 Puisque $\Psi$ (resp. $\Phi$) est asymptotique à $\widehat{Z}x^{B_{0}}$ (resp. $\widehat{Y}x^{A_{0}}$) dans un demi-plan droit, il existe des constantes $C, \ M >1$ telles que, pour tout $x\in\complex$ avec $Re(x) \geq M$, on ait :
\begin{equation*}
 \|(\Phi x^{-A_{0}}-\widehat{Y}_{|2N})x^{A_{0}} \omega x^{-B_{0}}+
\widehat{Y}_{|2N}x^{A_{0}} \omega x^{-B_{0}} -\widehat{Z}_{|2N}\| \leq
C |x|^{-2N-1}
\end{equation*}
et $$\| \Phi x^{-A_{0}}-\widehat{Y}_{|2N}\| \leq C |x|^{-2N-1}.$$
Il en résulte que, pour tout $x\in\complex$ avec $Re(x) \geq M$, on a :
\begin{equation*}
\| \widehat{Y}_{|2N}x^{A_{0}} \omega x^{-B_{0}} -\widehat{Z}_{|2N}\|
\leq C |x|^{-2N-1}+ C|x|^{-N-1} \|\omega\|.
\end{equation*}
Soit $S$ un segment compact
de $[M,+\infty[$ non réduit à un point, sur lequel $\omega$ n'a pas
de pôle; $\omega$ est donc bornée sur $S+\integer.$
On en déduit qu'il existe une constante $C'>0$ telle que, pour tout $x\in S+\integer$, on ait :
\begin{equation*}
\| \widehat{Y}_{|2N}x^{A_{0}} \omega x^{-B_{0}} -\widehat{Z}_{|2N}\| \leq  C'|x|^{-N-1}.
\end{equation*}
Par suite, compte tenu du fait que $\widehat{Y}_{|2N}$ est tangente à $I$ en
$+\infty,$  il existe une constante $C''>0$ telle que, pour tout $x\in S+\integer$, on ait :
\begin{equation*}
 \| \omega
-x^{-A_{0}}(\widehat{Y}_{|2N})^{-1}\widehat{Z}_{|2N}x^{B_{0}}\| \leq C'' |x|^{-1}.
\end{equation*}
Ceci implique que $\omega$ est constante.
\end{proof}

\begin{coro}\label{dd}
Soit $\Psi$ une solution du système (\ref{aa}) admettant un
déve\-loppement asymptotique dans un demi-plan droit de la forme
$\widehat{Z}x^{B_{0}}$ avec $\widehat{Z} \in
Gl_{n}(\complex[[x^{-1}]])$ tangent à $I$ en l'infini et $B_{0}
\in M_{n}(\complex).$ Alors, $B_{0}=A_{0}$ et $\Psi = e_{A}^{+}.$
\end{coro}

\begin{proof}
Nous savons que $\Psi=e^{+}_{A}N$ avec $N \in Gl_{n}(\complex).$
Notons $\widehat{Y}x^{A_{0}}$ le développement asymptotique de
$e_{A}^{+}$ (qui est aussi l'unique solution formelle de notre
système de la forme $\widehat{F}x^{K}$ avec $\widehat{F} \in
Gl_{n}(\complex[[x^{-1}]])$ tangent à $I$ en l'infini et $K \in
M_{n}(\complex)$) alors,
$\widehat{Y}x^{A_{0}}N=\widehat{Z}x^{B_{0}}$ donc
$\widehat{Z}x^{B_{0}}$ est solution formelle de notre système de
la forme (série formelle des $x^{-1}$ tangente à $I$ en
l'infini)$\cdot$(puissance de $x$) donc égale à
$\widehat{Y}x^{A_{0}},$ \textit{i.e.} $N=I.$
\end{proof}

\subsection{``Forme normale'' des solutions d'un sys\-tème singulier
régulier.}

\subsubsection{``Forme normale''.}

Donnons nous un système fuchsien quelconque :
\begin{equation} \label{aba}
\delta_{-1}Y=AY.
\end{equation}

En parfaite analogie avec les théories des équations
différentielles et aux \textit{q}-différences, la démarche pour
résoudre un système fuchsien général consiste à se ramener au cas
non résonnant via une transformation de jauge rationnelle. Le résultat suivant est prouvé dans \cite{harris bis}.

\begin{lem}
Soit $A \in M_{n}(\seriesfact)$. Soient $c_{1},...,c_{n}$ les
valeurs propres de $A(+\infty) \in M_{n}(\complex).$ Il existe une matrice $T$, à coefficients rationnels, de la forme $T=T_{0}+T_{1}x^{-1} \in
Gl_{n}(\complex(x))$ avec $T_{0}, \ T_{1} \in M_{n}(\complex)$, telle
que $A^T$ soit dans $M_{n}(\seriesfact)$ et que $A^{T}(+\infty)$ ait pour valeurs propres
$c_{1}+1,...,c_{n}.$
\end{lem}

Nous pouvons donc trouver $T \in Gl_{n}(\complex(x))$ tel que
$A^{T}$ soit holomorphe en l'infini et telle que le spectre de
$A^{T}(\infty)$ soit contenu dans $\mathcal{B}=\{0<Re(x)\leq 1\}$
($A^{T}$ est en particulier non résonnant). Il est ainsi possible de trouver un système fondamental de solutions, en $+\infty$, du
système fuchsien (\ref{aba}), de la forme
$M^{(+\infty)}N^{(+\infty)}$ où $M^{(+\infty)}=TF$ avec $F \in
Gl_{n}(\seriesfact)$ tangente à $I$ en $+\infty$ et
$N^{(+\infty)}$ est diagonale par blocs, ces derniers étant de la
forme :
\begin{equation*}
e_{cI_{\mu}+N_\mu}^{+}=\left( \begin{smallmatrix} l^{(0)}_{c}&l^{(1)}_{c}&...&l^{(\mu-1)}_{c}\\
                 \vdots&\ddots&&\vdots\\
                 \vdots&&l^{(0)}_{c}&l^{(1)}_{c}\\
                 0&...&0&l^{(0)}_{c}
\end{smallmatrix} \right), \ c\in \mathcal{B}.
\end{equation*}
On peut en outre faire en sorte que, pour un ordre fixé sur
$\complex$ (on peut par exemple prendre l'ordre total $\prec$
défini par $x \prec y \Longleftrightarrow [ Re(x)<Re(y) \text{ ou
}$ $ (Re(x)=Re(y) \text{ et } Im(x)<Im(y)) ] $), les exposants
soient rangés par taille croissante et que pour chaque exposant
$c$ fixé les blocs de Jordan soient également organisés de manière
croissante suivant leurs tailles.\\
\begin{rem}
Il suit d'une remarque antérieure
que $M^{(+\infty)}$ peut être choisi de telle sorte qu'elle soit, ainsi que son inverse, à croissance au plus
polynomiale dans tout demi-plan droit.
\end{rem}

Les conclusions précédentes subsistent clairement
pour les systèmes singuliers réguliers. Cela nous conduit à la définition suivante.

\begin{defi}
Une solution en $+\infty$ d'un système singulier régulier sera
dite sous forme normale si elle s'écrit comme un produit \linebreak
$M^{(+\infty)}N^{(+\infty)}$ avec $M^{(+\infty)} \in
Gl_{n}(\corpsseriesfact (\complex))$ et $N^{(+\infty)}$ diagonale
par blocs, de blocs diagonaux de la forme $e_{cI_{\mu}+N_\mu}^{+}$
avec $c \in \mathcal{B}$ et avec les $c$ rangés par ordre
croissant (pour l'ordre précédemment introduit par exemple) ainsi
que la taille des blocs de Jordan pour chaque $c$ fixé. La matrice
$N^{(+\infty)}$ est alors dite ``log-car'' sous forme normale.
\end{defi}

Les raisonnements ci-dessus montrent qu'il existe toujours une
solution sous forme normale. Il n'y a pas d'unicité pour les formes normales.

\begin{theo}\label{theo unicite forme normale}
Soit $Y=MN$ et $Y'=M'N'$ deux solutions en $+\infty$ sous forme
normale d'un même système singulier régulier. Alors $N'=N$ et il
existe $R \in Gl_{n}(\complex)$ qui commute avec $N$ telle que
$M'=MR.$
\end{theo}

\begin{proof}
Cette preuve reprend une $q$-analogue; voir \cite{sauloy  systemes
aux q diff singu reg }. Soit $R = M^{-1} M'$; c'est un élément de $Gl_{n}(\corpsseriesfact(\complex)).$ Il est clair que $N^{-1}RN'$
est 1-périodique, méromorphe sur $\complex.$ Notons
$e_{c_{i}}^{+}L_{i}$ les blocs diagonaux de $N;$ $L_{i}$ a ses
coefficients dans
$\complex[\widetilde{l}^{(0)}_{c_{i}},\widetilde{l}^{(1)}_{c_{i}},...]$
avec
$\widetilde{l}^{(k)}_{c_{1}}=\frac{l^{(k)}_{c_{i}}}{e^{+}_{c_{i}}}.$
Notations analogues pour $N'.$ On note $R_{i,j}$ les
blocs de $R$ compatibles avec les blocs de $N$ et de $N'.$ Les blocs de
$N^{-1}RN'$ sont les
$\frac{e^{+}_{c'_{j}}}{e^{+}_{c_{i}}}L_{i}^{-1}R_{i,j}L_{j}'$;
ils sont 1-périodiques. Ainsi, les coefficients de
$L_{i}^{-1}R_{i,j}L_{j}'$ sont des solutions de l'équation :
\begin{equation*}
\tau^{-1}y=\frac{x-1-c_{i}}{x-1-c'_{j}}y
\end{equation*}
à coefficients dans
$\corpsseriesfact(\complex)(\{\widetilde{l}^{(k)}_{c_{i}}\}_{k\geq
0 },\{\widetilde{l}^{(k)}_{c'_{j}}\}_{k\geq 0 }).$ Soit
$R_{i,j}=0$ soit $R_{i,j}\neq 0.$ Dans ce dernier cas, si on
note que, d'après les résultats de l'appendice,
$\widetilde{l}^{(k)}_{c_{i}},\ \widetilde{l}^{(k)}_{c'_{j}} \in
\corpsseriesfact(log(x)),$ alors on obtient $c_{i}=c'_{j}$
($L_{i}^{-1}R_{i,j}L_{j}' \in
Gl_{n}(\complex)$).

 On note $S_{i,j}=L_{i}^{-1}R_{i,j}L_{j}'$ qui est donc à
coefficients constants. En développant l'égalité
$L_{i}S_{i,j}=R_{i,j}L_{j}'$ dans la base (sur le corps
$\corpsseriesfact(\complex),$ \textit{cf.} appendice) des
$\widetilde{l}^{(k)}_{c'_{i}},$ on déduit (si on note \linebreak
$L_{i}=\sum_{k\geq 0}\widetilde{l}^{(k)}_{c_{i}}\xi^{k}_{0,m_{i}}$
et $L'_{i}=\sum_{k\geq
0}\widetilde{l}^{(k)}_{c'_{i}}\xi^{k}_{0,m'_{i}}$) que
$R_{i,j}=S_{i,j}$ et $\xi^{k}_{0,m_{i}}R_{i,j}$ $=$
$R_{i,j}\xi^{k}_{0,m'_{i}}.$ Ainsi, $R \in Gl_{n}(\complex)$ et
$R^{-1}KR=K'$ avec $K=(\tau^{-1}N)N^{-1}$ et
$K'=(\tau^{-1}N')N'^{-1}.$ Les conditions de normalisation
impliquent alors $K=K'.$
\end{proof}
Le résultat suivant découle de remarques antérieures.
\begin{prop}\label{cc}
Si $M^{(+\infty)}N^{(+\infty)}$ est une solution canonique sous
forme normale alors $M^{(+\infty)}$ est à croissance au plus
polynomiale dans tout demi-plan droit. De même pour
$(M^{(+\infty)})^{-1}.$
\end{prop}

On a bien sûr des résultats similaires en $-\infty$.

\subsubsection{Caractérisation des systèmes singuliers réguliers
par la forme des solutions.}

\begin{lem} \label{lemme caract des syst par les sols}
Toute matrice $M \in Gl_{n}(\corpsseriesfact)$ peut s'écrire sous
la forme $M=CR$ avec $C \in Gl_{n}(\complex(x))$ et $R \in
Gl_{n}(\seriesfact)$ tangent à $I$ en $+\infty.$
\end{lem}

\begin{proof}
La preuve est laissée au lecteur.
Elle est une adaptation facile de celle du résultat analogue pour
les $q$-différences \cite{sauloy systemes aux q diff singu
reg }.
\end{proof}

\begin{theo}\label{caracterisation par la forme des solutions}
Si un système algébrique $\tau_{-1}Y=AY$ admet une solution $MN$
avec $M \in Gl_{n}(\corpsseriesfact)$ et $N$ une matrice ``log-car''
alors il est singulier régulier.
\end{theo}

\begin{proof}
On écrit $M=CR$ comme dans le lemme précédent et on définit
$B=(\tau_{-1}C)^{-1}AC \in M_{n}(\complex(x)).$ Le système défini
par $A$ est rationnellement équivalent à celui défini par $B.$ Mais
$B = (\tau_{-1}R) (\tau_{-1}NN^{-1}) R^{-1}.$ Ainsi $B(\infty)=I$ et
$B$ défini un système fuchsien.
\end{proof}

\section{Matrice de connexion de Birkhoff et classification rationnelle.}

\subsection{Notations}
Soit (en notant $\mathcal{M}(\complex / \rinteger)$ le corps des
fonctions méromorphes sur $\complex$ et 1-périodiques) :
\begin{equation*}
  \mathcal{A}_{n}'=Gl_{n}(\mathcal{M}(\complex / \rinteger))\times\{ \text{matrices ``log-car'' sous forme normale}\}.
\end{equation*}
On définit la relation d'équivalence $\sim$ sur
$\mathcal{A}_{n}'$ par :
\begin{equation*}
(P,N) \sim (P',N')
\end{equation*}
\begin{equation*}
\Leftrightarrow
\end{equation*}
\begin{equation*}
 N=N' \text{ et }
\end{equation*}
\begin{equation*}
 \exists R_{1}, \ R_{2} \in
Gl_{n}(\complex) \text{ t.q. } R_{1}P=P'R_{2} \text{ et }
[R_{i},N]=0, \ i=1,2;
\end{equation*}
la classe d'un couple $(P,N)$ est notée $\overline{(P,N)}.$ On
pose :
\begin{equation*}
  \mathcal{F}_{n}'=\mathcal{A}_{n}'/\sim.
\end{equation*}
L'ensemble suivant jouera un rôle fondamental :
\begin{equation*}
  \mathcal{F}_{n}=\{\overline{(P,N)} \ |\ P\in Gl_{n}(\complex(e^{2\pi i x}))\text{ et } P(\pm i \infty)=e^{\pm i \pi N_{0} }\} \subset \mathcal{F}_{n}'
\end{equation*}
où l'on note(ra) $N_{0} \in M_{n}(\complex)$ la réduite de Jordan
sous forme normale correspondant à $N$ \textit{i.e.} la réduite de
Jordan $N_{0}$ telle que $e^{+}_{N_{0}}=N.$

Enfin, l'ensemble des classes de systèmes singuliers réguliers
modulo la relation d'équivalence rationnelle est noté
$\mathcal{E}_{n}.$ On rappelle que deux systèmes singuliers
réguliers $\tau_{-1} Y=AY$ et $\tau_{-1} Y=BY$ sont dits rationnellement
équivalents s'il existe une transformation de jauge rationnelle
qui transforme le premier en le second \textit{i.e.} s'il existe $R\in Gl_n(\complex(x))$ telle que $B(x)=R(x-1)^{-1}A(x)R(x)$; c'est clairement
une relation d'équivalence.

Notre but est de décrire les classes de systèmes pour
l'équivalence rationnelle grâce à la matrice de connexion de
Birkhoff. Le principe remonte à G.D. Birkhoff.

\subsection{Définition et propriétés.}

Considérons un système singulier régulier. Soient :
\begin{equation*}
Y^{(-\infty)}=M^{(-\infty)}N^{(-\infty)} \text{ et }
Y^{(+\infty)}=M^{(+\infty)}N^{(+\infty)}
\end{equation*}
des solutions sous forme normale, en $-\infty$ et en $+\infty$
respectivement. Soit :
\begin{eqnarray*}
P&=&(Y^{(+\infty)})^{-1}Y^{(-\infty)}\\
&=&(N^{(+\infty)})^{-1}(M^{(+\infty)})^{-1}M^{(-\infty)}N^{(-\infty)} \in Gl_{n}(\mathcal{M}(\complex /
\rinteger))
\end{eqnarray*}
la matrice de connexion associée. D'après le théorème \ref{theo unicite forme
normale}, on définit ainsi une application :
\begin{equation*}
  Bir : \mathcal{E}_{n} \longrightarrow \mathcal{F}_{n}';
\end{equation*}
elle associe à la classe d'une équation, la classe du couple formé
d'une
matrice de connexion et de la réduite de Jordan sous forme normale correspondant à $N^{(+\infty)}.$

En outre, d'après la proposition \ref{cc}, $M^{(-\infty)}$ et
$(M^{(+\infty)})^{-1}$ sont à croissance au plus polynomiale dans
les bandes verticales. Il en va de même pour $N^{(-\infty)}$ et
$(N^{(+\infty)})^{-1}$ (d'après des propriétés classiques de la
fonction Gamma). Ainsi $(Y^{(+\infty)})^{-1}Y^{(-\infty)}$ est une
fonction méro\-morphe 1-périodique à croissance modérée dans les
bandes verticales et donc une fonction trigonométrique.

Enfin, rappelons qu'on peut trouver deux solutions sous forme
normale qui s'écrivent :
\begin{equation*}
Y^{(-\infty)}=Te^{-}_{A} \text{ et }
Y^{(+\infty)}=Te^{+}_{A}
\end{equation*}
où $A$ est holomorphe à l'infini et $A(\infty)=N_{0}$ est la
matrice de Jordan correspondant à la forme normale de notre
système et où $T \in Gl_{n}(\complex(x))$ (une transformation de
jauge rationnelle qui ramène au cas fuchsien non résonnant).
Ainsi, $(Y^{(+\infty)})^{-1}Y^{(-\infty)}=(e^{+}_{A})^{-1}e^{-}_{A}$
et les propriétés asymptotiques de $e^{-}_{A}$ et de $e^{+}_{A}$
données lors de la sous-section \ref{sous section comportement
asymptotique} impliquent :
\begin{equation*}
(Y^{(+\infty)})^{-1}Y^{(-\infty)}(\pm i \infty)=e^{\pm i \pi
N_{0}}.
\end{equation*}
Finalement, on peut co-restreindre le but de l'application
$Bir$ à $\mathcal{F}_{n}$. On
note encore $Bir$ l'application obtenue.\\

\subsection{Matrice de Birkhoff et classification rationnelle.}

Le théorème principal, dont le principe remonte à Birkhoff, peut
maintenant être énoncé.

\begin{theo}
L'application $Bir : \mathcal{E}_{n} \longrightarrow \mathcal{F}_{n}$
est une bijection.
\end{theo}

\begin{proof}

\textsl{-Injectivité.}

Soient deux systèmes $\tau_{-1}Y=AY$ et  $\tau_{-1}Y=BY$ donnant le même
élément de $\mathcal{F}_{n}$ par $Bir.$ On peut donc choisir :
\begin{equation*}
Y^{(-\infty)}_{A,B}=M^{(-\infty)}_{A,B}N^{(-\infty)}_{A,B}
\text{ et }
Y^{(+\infty)}_{A,B}=M^{(+\infty)}_{A,B}N^{(+\infty)}_{A,B}
\end{equation*}
des solutions canoniques sous forme normale en $-\infty$ et en
$+\infty$ respectivement telles que
$(Y^{(+\infty)}_{A})^{-1}Y^{(-\infty)}_{A}=(Y^{(+\infty)}_{B})^{-1}Y^{(-\infty)}_{B}$.
Par suite :
\begin{equation*}
  Y^{(+\infty)}_{A}(Y^{(+\infty)}_{B})^{-1}=Y^{(-\infty)}_{A}(Y^{(-\infty)}_{B})^{-1}.
\end{equation*}
Rappelons que $N^{(\pm \infty)}_{A}=N^{(\pm \infty)}_{B}.$ D'une
part, $Y^{(+\infty)}_{A}(Y^{(+\infty)}_{B})^{-1}$ est à croissance
au plus polynomiale dans tout demi-plan droit (cela résulte de
remarques précédentes), d'autre part
$Y^{(+\infty)}_{A}(Y^{(+\infty)}_{B})^{-1}$ est à croissance
au plus polynomiale dans tout demi-plan gauche par l'égalité
ci-dessus; ainsi cette matrice, \textit{a priori} seulement
méromorphe sur $\complex,$ est en réalité rationnelle; puisqu'elle
conjugue nos deux systèmes, cela prouve l'injectivité de
$Bir.$\\

\textsl{-Surjectivité.}

Il suffit de reprendre la méthode de G.D. Birkhoff dans
\cite{Birkhoff generalised RH pb}. Voici quelques indications, en
reprenant ses notations. Soient $(P,N)$ un
couple de $\mathcal{F}_{n}$ et $N_{0}$ la matrice sous forme normale
correspondant à $N;$ on commence comme Birkhoff (paragraphe 17) : $A_{1}(x)=TPT^{-1}$ où $T(x)=x^{N_{0}}$ avec
les mêmes choix de branches du logarithme. Alors $A_{1}$ est asymptotique à
$I$ dans tout secteur d'angle $< \pi$ bisecté par le demi-axe des
imaginaires purs supérieur; de même en bisectant par le demi-axe
des imaginaires purs inférieur. Nous pouvons maintenant utiliser
les mêmes raisonnements que G.D. Birkhoff : les deux mêmes
utilisations de son ``preliminary theorem''
permettent de prouver l'existence de deux fonctions méromorphes
$Y^{-}$ et $Y^{+}$ telles que $Y^{-}=Y^{+}P$ et qui sont
asymptotiques dans un demi-plan droit pour l'une, gauche pour
l'autre, à une fonction $S$ de la forme $S=Yx^{N_{0}}$ où $Y \in
Gl_{n}(\complex((x^{-1}))).$

\begin{lem}
Toute matrice $M \in Gl_{n}(\complex((x^{-1})))$ peut s'écrire sous la
forme $M=CR$ avec $C \in Gl_{n}(\complex(x))$ et $R \in
Gl_{n}(\complex[[x^{-1}]])$ tangent à $I$ en $\infty.$
\end{lem}

\begin{proof}
C'est essentiellement la même que celle du lemme \ref{lemme caract
des syst par les sols}.
\end{proof}

D'après le lemme, il existe deux fonctions méromorphes $Y^{-}$ et
$Y^{+}$ telles que $Y^{-}=Y^{+}P$ et qui sont asymptotiques, dans
un demi-plan droit pour l'une, gauche pour l'autre, à $S$ de la
forme $S=Yx^{A_{0}}$ où $Y \in Gl_{n}(\complex[[x^{-1}]])$ est tangent
à $I$ en l'infini. Comme Birkhoff, on en conclut que
$Q=Y^{-}(x)(Y^{-}(x-1))^{-1}=Y^{-}(x)(Y^{-}(x-1))^{-1}$ est
rationnelle. Pour terminer remarquons que $Q$ est de la forme
$I+A_{0}/x+(\text{termes de degrés plus petits}).$ Ainsi, on a
trouvé un système (rationnel) fuchsien, défini par une matrice $A$
telle que $A(\infty)=N_{0}$ (le système est donc non résonnant) et
admettant deux solutions $Y^{-}$ et $Y^{+}$ telles que
$Y^{-}=Y^{+}P$ et qui sont asymptotiques, dans un demi-plan droit
pour l'une, gauche pour l'autre, à $S$ de la forme $S=Yx^{A_{0}}$
où $Y \in Gl_{n}(\complex[[x^{-1}]])$ est tangent à $I$ en l'infini. Le
corollaire \ref{dd} implique que $Y^{-}$ et $Y^{+}$ sont deux
solutions canoniques. Ceci termine la preuve de la surjectivité de
l'application $Bir.$
\end{proof}

\section{Confluence des systèmes aux diffé\-rences fuchsiens non réson\-nants vers les sys\-tèmes différentiels.}

Pour les notations et les terminologies employées,
nous renvoyons le lecteur à la section \ref{notations et
terminologies}. \\

Nous nous
intéres\-sons à partir de maintenant, et jusqu'à la fin de cet
article, à la confluence des systèmes aux différences fuchsiens
non résonnants vers les systèmes différentiels. Autrement dit,
nous étudions les propriétés d'un système différentiel fuchsien en
l'infini, en l'occurrence :
\begin{equation*}
\delta \widetilde{Y}=\widetilde{A}\widetilde{Y};
\end{equation*}
($\delta$ désignant l'opérateur d'Euler $\delta = t \frac{d}{dt}$)
que l'on peut (re)trouver en le voyant comme limite quand $h
\longrightarrow 0$ de systèmes aux différences de pas $h>0$ :
\begin{equation*}
\delta_{-h}Y^{(h)}=A^{(h)}Y^{(h)}.
\end{equation*}
Nous considérons des systèmes de pas $h$ plus
généraux que les systèmes algébriques : on les supposera définis par
des fonctions développables en séries de $h$-factorielles
(ou de $h$-rétro-factorielles suivant qu'on s'intéresse
à $\pm \infty$).

\subsection{Résolution des systèmes à matrice cons\-tante.}
\label{resolution des systeme a matrice constante} On généralise
facilement notre résolution des équations des caractères et des
logarithmes à des équations de pas $h$ quelconque.

Pout tout $c\in \complex$ et tout $h>0$ nous nous donnons une
famille $l^{(k,h)}_{c}, \ k\in \integer$, de fonctions uniformes
et méromorphes sur $\complex$ satisfaisant :
\begin{equation*}
\delta_{-h}l^{(k,h)}_{c}=cl^{(k,h)}_{c}+l^{(k-1,h)}_{c}.
\end{equation*}
Posons $e^{(h)+}_{c}=l^{(0,h)}_{c}$ et $l^{(k,h)}=l^{(k,h)}_{0}.$
Ces fonctions vérifient les équations suivantes :
\begin{eqnarray*}
\delta_{-h}e^{(h)+}_{c}&=&ce^{(h)+}_{c}  \text{ (équations des caractères),}\\
\delta_{-h}l^{(k,h)}&=&l^{(k-1,h)} \text{ (équations des logarithmes).}
\end{eqnarray*}

Il faut choisir des solutions qui présentent de bonnes propriétés
de confluence. Voici des choix possibles.\\
\\
\begin{ex}\textbf{\ 1.}
Soient :
\begin{equation*}
e_{c}^+(x)=\frac{\Gamma(x)}{\Gamma(x-c)},\ \
e^{(h)+}_{c}(x)=h^{c}e_{c}^+(\frac{x}{h}),\ \
l^{(k,h)}_{c}(x)=\frac{1}{k!}\frac{\partial^{k}}{\partial
c'^{k}}_{|c'=c} e^{(h)+}_{c'}(x).
\end{equation*}

\begin{prop}
 Soient $x \in \complex \backslash \real^{-}$ et, pour $h\in \io$, $c^{(h)} \in \complex$ tels que $c^{(h)} \xrightarrow[h \longrightarrow 0^{+}]{}
c\in \complex$.
Alors, pour tout $k\in \integer$:
\begin{equation*}
l^{(k,h)}_{c^{(h)}}(x) \xrightarrow[h \longrightarrow 0^{+}]{}
x^{c}\frac{\log^{k}(x)}{k!}.
\end{equation*}
\end{prop}
\begin{proof}
\begin{lem}\label{lemme conv caract}
Soit $x \in \complex \backslash \real^{-}.$ Alors :
\begin{equation*}
  e^{(h)+}_{c}(x) \xrightarrow[h \longrightarrow 0^{+}]{}  x^{c}
\end{equation*}
uniformément en $c \in D(0,r)=\{y \in \complex \ | \ |y|<r  \}.$
\end{lem}
\begin{proof}
Résulte facilement de la formule de Stirling.
\end{proof}
Fixons $x\in \complex \backslash \real^{-}$ et $r>0$. On
considère la famille (indexée par $h$) de fonctions holomorphes sur
$\overline{D}(0,r)$ $f_{h} : c \longmapsto e^{(h)+}_{c}(x)$. D'après le lemme \ref{lemme conv caract}, la suite de fonctions $f_{h}$ converge uniformément vers $f : c \longmapsto
x^{c}$ lorsque $h \rightarrow 0^+$. Le théorème de Weierstrass
implique que, pour tout $k\in \integer$, $f_{h}^{(k)}$
converge uniformément vers $f^{(k)}$. Nous en déduisons que, pour tout $k\in \integer$,
$f_{h}^{(k)}(c^{(h)}) \xrightarrow[h \longrightarrow 0^{+}]{}
f^{(k)}(c)$. Ce qui est le résultat attendu.
\end{proof}
\end{ex}
$_{}$\\
\begin{ex}\textbf{\ 2.}
Soient :
\begin{equation*}
e_{c}^+(x)=\frac{\Gamma(1+c-x)}{\Gamma(1-x)},\ \
e^{(h)+}_{c}(x)=h^{c}e_{c}^+(\frac{x}{h}),\ \
l^{(k,h)}_{c}(x) = \frac{1}{k!}\frac{\partial^{k}}{\partial
c'^{k}}_{|c'=c} e^{(h)+}_{c'}(x).
\end{equation*}

Comme dans l'exemple précédent, on prouve la
\begin{prop}
 Soit $x \in \complex \backslash \real^{+}$ et, pour $h\in \io$, $c^{(h)} \in \complex$ tels que $c^{(h)} \xrightarrow[h \longrightarrow 0^{+}]{}
c\in \complex$.
Alors, pour tout $k\in \integer$ :
\begin{equation*}
l^{(k,h)}_{c^{(h)}}(x) \xrightarrow[h \longrightarrow 0^{+}]{}
(-x)^{c}\frac{\log^{k}(-x)}{k!}.
\end{equation*}
\end{prop}

\end{ex}

Pour résoudre un système constant, on procède comme dans le cas du
pas $h=1.$ On se donne donc un système :
\begin{equation*}
\delta_{-h}Y=A^{(h)}_{0}Y
\end{equation*}
où $A^{(h)}_{0} \in M_{n}(\complex)$. On réduit $A^{(h)}_0$ sous forme de
Jordan :
\begin{equation*}
A^{(h)}_{0}=P^{(h)} diag(c_{1}I_{\mu_{1}}+N_{\mu_1},...,c_{m}I_{\mu_{m}}+N_{\mu_m})\left(P^{(h)}\right)^{-1}
\end{equation*}
suivant la même règle que dans le cas du pas $h=1$.
Alors, un système fondamental de solutions, en $+\infty$, du système de pas $h$
défini par $A^{(h)}_{0}$ est donné par:

\begin{equation*}
e_{A_{0}^{(h)}}^{(h)+}:=P^{(h)} diag(e_{c_{1}I_{\mu_{1}}+N_{\mu_1}}^{(h)+},...,e_{c_{m}I_{\mu_{m}}+N_{\mu_m}}^{(h)+})\left(P^{(h)}\right)^{-1}
\end{equation*}
avec :
\begin{equation*}
e_{cI_{\mu}+N_\mu}^{(h)+}=\left( \begin{smallmatrix} l^{(0,h)}_{c}&l^{(1,h)}_{c}&...&l^{(\mu-1,h)}_{c}\\
                 \vdots&\ddots&&\vdots\\
                 \vdots&&l^{(0,h)}_{c}&l^{(1,h)}_{c}\\
                 0&...&0&l^{(0,h)}_{c}
\end{smallmatrix} \right).
\end{equation*}
La
solution obtenue est indépendante de la réduction de Jordan
choisie. Nous aurons besoin d'une hypothèse sur les
réductions de Jordan (analogue à une hypothèse de \cite{sauloy systemes aux q diff singu
reg }).

\begin{defi} \label{def deploiement}
Soient $\widetilde{A}_{0} \in M_{n}(\complex)$ et, pour $h\in\io$,
$A_{0}^{(h)} \in M_{n}(\complex)$. On dira qu'une réduction de
Jordan de $\widetilde{A}_{0}$ se déploie en des réductions de
Jordan des $A_{0}^{(h)},\ h>0,$ s'il existe une réduction de
Jordan de $\widetilde{A}_{0}$ :
\begin{equation*}
\widetilde{A}_{0}=\widetilde{P}\widetilde{J}\widetilde{P}^{-1}
\end{equation*}
et des réductions de Jordan des $A_{0}^{(h)},$ $h\in\io$ :
\begin{equation*}
A_{0}^{(h)}=P^{(h)}J^{(h)}\left(P^{(h)}\right)^{-1}
\end{equation*}
telles que :
\begin{equation*}
P^{(h)} \xrightarrow[h\longrightarrow 0^+]{}\widetilde{P},\ J^{(h)}
\xrightarrow[h\longrightarrow 0^+]{}\widetilde{J}.
\end{equation*}
\end{defi}
Notons que l'hypothèse de déploiement implique que $A_{0}^{(h)}
\xrightarrow[h\longrightarrow 0^+]{}\widetilde{A}_{0}.$ Compte tenu
des résultats précédents, on a clairement le résultat suivant.

\begin{prop} \label{conflu dans le cas des coeff constants}
Soit $\widetilde{A}_{0} \in M_{n}(\complex)$ et, pour $h\in\io$,
$A_{0}^{(h)} \in M_{n}(\complex)$. On suppose qu'une réduction de
Jordan de $\widetilde{A}_{0}$ se déploie en des réductions de
Jordan des $A_{0}^{(h)},\ h\in\io.$ Alors, dans le cas de l'exemple 1
(resp. 2) on a, $\forall x \in \complex \setminus  \real^-$ (resp.
$\complex \setminus \real^+$) :
$$e_{A_{0}^{(h)}}^{(h)+}(x) \xrightarrow[h\longrightarrow 0^+]{} x^{\widetilde{A}_{0}} \text{ (resp. $(-x)^{\widetilde{A}_0}$)}.$$
\end{prop}

\subsection{On se ramène au cas constant.}
\subsubsection{La transformation de jauge: énoncé du théorème.}\label{qsd}

Soit $\|\cdot\|$ une norme d'algèbre sur $M_{n}(\complex).$
On note $|||\cdot|||$ la norme (sur l'espace vectoriel des
endomorphismes de $M_{n}(\complex)$) subordonnée à $\|\cdot\|$.
Désignons enfin, pour $s \in \integer^*$, par $K_{s,A_{0}^{(h)}}$ l'opérateur
linéaire sur $M_{n}(\complex)$ défini par :
$$\forall M \in M_{n}(\complex),
K_{s,A_{0}^{(h)}}(M)=A^{(h)}_{0}M-MA^{(h)}_{0}-sM.$$

Cette partie est consacrée à la preuve du théorème suivant.

\begin{theo}
Soit $\delta_{-h}Y=A^{(h)}Y$, $h\in\io$ une famille de
systèmes de Fuchs $(C,\lambda)$ en $+\infty,$  non résonnants
(non nécessairement algébriques). On suppose de plus que :
$$\sup_{s\in \integer^{*}, \ h\in \io}|||K_{s,A_{0}^{(h)}}^{-1}|||<\infty.$$ Alors, $\forall h\in\io,\ \exists !  F^{(h)} \in
Gl_n(\seriesfact^{(h)})$ tel que $F^{(h)}(+\infty)=I$ et
$(A^{(h)})^{F^{(h)}}=A^{(h)}_{0}.$ De plus, la famille des
$F^{(h)}$, $h\in\io$ est d'un certain type $(C',\lambda').$
\end{theo}

La preuve que nous donnons de ce théorème est inspirée de techniques utilisées par A. Duval dans
\cite{duvseriesfact} et \cite{duvconfqfuchsien}.\\

\textbf{Notations. }Si $A^{(h)}(x)=\sum_{s=0}^{+\infty}
A^{(h)}_{s}x^{-[s]_{h}} \in M_n(\seriesfact^{(h)})$, on note :
\begin{equation*}
\overline{A}^{(h)}(x)=A^{(h)}(hx)=\sum_{s=0}^{+\infty}
\overline{A}^{(h)}_{s}x^{-[s]}
\end{equation*}
avec :
\begin{equation*}
\overline{A}^{(h)}_{s}=\frac{A^{(h)}_{s}}{h^{s}}.
\end{equation*}

\subsubsection{Préparatifs.}

\begin{prop}\label{preparatif}
Soit $\delta_{-h}{Y}=A^{(h)}Y$, $h\in\io$ une famille de
systèmes de Fuchs $(C,\lambda)$ en $+\infty$ (non nécessairement algé\-bri\-ques), tels que, pour $h\in\io$, le spectre de $A^{(h)}_0=A^{(h)}(+\infty)$ ne contienne aucun entier strictement négatif. Soient $b=\sup_{s\in
\integer^{*},h\in\io} \|(sI+A^{(h)}_{0})^{-1}\|$ qu'on suppose fini
et $0 \neq U_{0} \in \complex^{n}.$ Alors, pour tout $h\in\io,$ l'équation $\delta_{-h}{Y}=A^{(h)}Y$ admet une
solution dans $M_{n,1}(\seriesfactform^{(h)})$ de terme constant
$U_{0}$ si et seulement si $U_{0}$ est dans le noyau de
$A^{(h)}_{0}.$ Lorsqu'elle existe, une telle
solution est unique, converge dans un demi-plan droit, et la
famille de ces solutions est de type $(Cb\|U_{0}\|,Cb+\lambda).$
\end{prop}

\begin{proof} Notons :
\begin{equation*}
A^{(h)}(x)=\sum_{s=0}^{+\infty} A^{(h)}_{s}x^{-[s]_{h}}
\text{ et }
Y^{(h)}(x)=\sum_{s=0}^{+\infty} Y^{(h)}_{s}x^{-[s]_{h}}.
\end{equation*}

Rappelons l'hypothèse :

\begin{equation*}
\|\overline{A}^{(h)}_{s}\| \leq \frac{C}{h}
\left(\frac{\lambda}{h}\right)^{[s-1]} =
\overline{C}^{(h)}\left(\overline{\lambda}^{(h)}\right)^{[s-1]},\
s \in \integer^{*} 
\end{equation*}
avec :
\begin{equation*}
\overline{C}^{(h)}=\frac{C}{h} \text{ et } \overline{\lambda}^{(h)}=\frac{\lambda}{h}.
\end{equation*}

\textsl{Partie formelle. } Tout revient à trouver $Y^{(h)}\in
M_{n,1}(\seriesfactform^{(h)})$ de ``terme constant'' $U_{0}$ et
telle que $\overline{Y}^{(h)}$ soit solution de
$\delta_{-1}{\overline{Y}^{(h)}}(x)=\overline{A}^{(h)}(x)\overline{Y}^{(h)}(z).$
Nous avons (formule de multiplication de deux séries de
factorielles):
\begin{equation*}
\overline{A}^{(h)}(x)\overline{Y}^{(h)}(x)=\sum_{s=0}^{+\infty}C^{(h)}_{s}x^{-[s]}
\end{equation*}
avec:
\begin{equation*}
C^{(h)}_{0}=\overline{A}^{(h)}_{0}\overline{Y}^{(h)}_{0}; \
C^{(h)}_{s}=\overline{A}^{(h)}_{0}\overline{Y}^{(h)}_{s}+\overline{A}^{(h)}_{s}\overline{Y}^{(h)}_{0}+\sum_{(j,k,l)
\in J_{s}
}c_{j,l}^{(k)}\overline{A}^{(h)}_{j}\overline{Y}^{(h)}_{l}
\end{equation*}
et:
\begin{equation*}
J_{s} = \{(j,k,l) \backslash j,l \geq 1,\ k \geq 0,\ j+k+l=s\}
\end{equation*}
\begin{equation*}
c_{j,l}^{(k)}=\frac{(j+k-1)!(l+k-1)!}{k!(j-1)!(l-1)!}.
\end{equation*}
D'autre part :
\begin{equation*}
\delta_{-1}\overline{Y}^{(h)}(x)=\sum_{s=0}^{+\infty} -s
\overline{Y}_{s}^{(h)}x^{-[s]}.
\end{equation*}
Ainsi, $\overline{Y}^{(h)}$ est solution si et seulement si :
\begin{equation*}
\begin{cases}
    \overline{A}^{(h)}_{0}\overline{Y}^{(h)}_{0}=0  \\
    -(sI+\overline{A}^{(h)}_{0})\overline{Y}^{(h)}_{s}=\overline{A}^{(h)}_{s}Y^{(h)}_{0}+\sum_{(j,k,l) \in
J_{s} }c_{j,l}^{(k)}\overline{A}^{(h)}_{j}\overline{Y}^{(h)}_{l},\
s \in \integer^{*}.
  \end{cases}
  \end{equation*}
Il est donc nécessaire que
$\overline{Y}^{(h)}_{0}$ soit dans $Ker(\overline{A}^{(h)}_{0}).$
Dans ce cas, le système se résout de
proche en proche :
\begin{equation*}
    \overline{Y}^{(h)}_{s}=-(sI+\overline{A}^{(h)}_{0})^{-1}
    \left[\overline{A}^{(h)}_{s}\overline{Y}^{(h)}_{0}+\sum_{(j,k,l) \in
J_{s} }c_{j,l}^{(k)}\overline{A}^{(h)}_{j}\overline{Y}^{(h)}_{l}
\right].
  \end{equation*}
Ceci clôt l'aspect formel de la proposition.\\

\textsl{Partie convergente. } Voyons maintenant ce qu'il en est de
la convergence de la solution obtenue.

Notons : $\overline{a}^{(h)}_{s}=\|\overline{A}^{(h)}_{s}\|$ et
$\overline{y}^{(h)}_{s}=\|\overline{Y}^{(h)}_{s}\|$, $\forall s
\in
\integer.$ \\
De la relation de récurrence ci-dessus, on déduit :
\begin{equation*}
    \overline{y}^{(h)}_{s}\leq b\left[\overline{a}^{(h)}_{s}\overline{y}^{(h)}_{0}
    +\sum_{(j,k,l) \in
J_{s} }c_{j,l}^{(k)}\overline{a}^{(h)}_{j}\overline{y}^{(h)}_{l}
\right].
  \end{equation*}
Introduisons une série majorante :
\begin{equation*}
  \begin{cases}
   \overline{y}^{(h)>}_{1}=b\overline{a}^{(h)}_{1}\overline{y}^{(h)}_{0} \\
    \overline{y}^{(h)>}_{s}= b\left[\overline{a}^{(h)}_{s}\overline{y}^{(h)}_{0}
    +\sum_{(j,k,l) \in
J_{s} }c_{j,l}^{(k)}\overline{a}^{(h)}_{j}\overline{y}^{(h)>}_{l}
\right],\ s \in \integer^{*}.
  \end{cases}
  \end{equation*}
Elle domine la série de terme général $\overline{y}^{(h)}_{s}.$
Notons:
\begin{equation*}
  \overline{y}^{(h)>}(x)=\sum_{s=1}^{+\infty}
\overline{y}^{(h)>}_{s}x^{-[s]} \in \seriesfactform
  \end{equation*}
 et:
 \begin{equation*}
  \overline{a}^{(h)>}(x)=\sum_{s=1}^{+\infty}
\overline{a}^{(h)>}_{s}x^{-[s]} \in \seriesfactform.
  \end{equation*}
La formule de multiplication de deux séries de factorielles montre
que:
\begin{equation*}
 \overline{y}^{(h)>}(x)=b(\overline{y}^{(h)}_{0}\overline{a}^{(h)}(x)
+\overline{a}^{(h)}(x)\overline{y}^{(h)>}(x)),
  \end{equation*}
par conséquent :
\begin{equation*}
 \overline{y}^{(h)>}(x)
=-\overline{y}^{(h)}_{0}\left(
1-\frac{1}{1-b\overline{a}^{(h)}(x)}\right).
  \end{equation*}
A présent, on utilise les estimations sur les coefficients de
l'inverse d'une série de factorielles d'un type $(C,\lambda)$
(voir \cite{duvseriesfact} ou \cite{Norlund series d interp});
elles donnent $\overline{y}^{(h)>}_{s} \leq
\overline{C}^{(h)}b\overline{y}^{(h)}_{0}
\frac{\Gamma(\overline{\lambda}^{(h)}+\overline{C}^{(h)}b+s-1)}{\Gamma(\overline{\lambda}^{(h)}+\overline{C}^{(h)}b)}.$
\end{proof}

\subsubsection{La transformation de jauge: preuve du théorème.  }

Rappelons l'hypothèse :

\begin{equation*}
\|\overline{A}^{(h)}_{s}\| \leq \frac{C}{h}
\left(\frac{\lambda}{h}\right)^{[s-1]} = \overline{C}^{(h)}
(\overline{\lambda}^{(h)})^{[s-1]}
\end{equation*}
avec :
\begin{equation*}
\overline{C}^{(h)}=\frac{C}{h} \text{ et } \overline{\lambda}^{(h)}=\frac{\lambda}{h}.
\end{equation*}

Passons à présent à la preuve du théorème. Notons
$F^{(h)}(x)=\sum_{s=0}^{+\infty} F^{(h)}_{s}x^{-[s]_{h}}.$ La
matrice $F^{(h)}$ satisfait les conditions de la proposition si et
seulement si $A^{(h)}(x)F^{(h)}(x)-\delta_{-h} F^{(h)}(x)=F^{(h)}(x-h)A_{0}^{(h)}.$

Ainsi, tout revient à trouver $F^{(h)} \in Gl_{n}(\seriesfact^{(h)})$
(tangente à $I$ en $+\infty$) telle que $\overline{F}^{(h)}$ soit
développable en série de factorielles tangente à l'identité en
l'infini et telle que
$\overline{A}^{(h)}(x)\overline{F}^{(h)}(x)-\delta_{-1}
\overline{F}^{(h)}(x)=\overline{F}^{(h)}(x-1)\overline{A}^{(h)}_{0}.$

Cette dernière équation peut se réécrire comme suit :

\begin{equation*}
  \delta_{-1} \overline{F}^{(h)}(x)=\sum_{s=0}^{+\infty}L^{(h)}_{s}(\overline{F}^{(h)}(x)) x^{-[s]}
\end{equation*}
où, pour $s=0,1$, $L^{(h)}_{s}$ est l'opérateur linéaire défini par :
\begin{equation*}
  L^{(h)}_{0}(U)=\overline{A}^{(h)}_{0}U-U\overline{A}^{(h)}_{0},
\end{equation*}
\begin{equation*}
  L^{(h)}_{1}(U)=\overline{A}^{(h)}_{1}(U)+(\overline{A}^{(h)}_{0}U-U\overline{A}^{(h)}_{0})\overline{A}^{(h)}_{0}
\end{equation*}
et pour $s \geq 2$ :
\begin{equation*}
  L^{(h)}_{s}(U)=\overline{A}^{(h)}_{s}(U)+(\overline{A}^{(h)}_{0}U-U\overline{A}^{(h)}_{0})B^{(h)}_{s}+\sum_{(j,k,l) \in
J_{s} }c_{j,l}^{(k)}\overline{A}^{(h)}_{j}UB^{(h)}_{l}
\end{equation*}
avec:
\begin{equation*}
  B^{(h)}_{s}=\overline{A}^{(h)}_{0}(\overline{A}^{(h)}_{0}+I)\cdots (\overline{A}^{(h)}_{0}+(s-1)I),\ s \geq 1.
\end{equation*}
Puisque $\|\overline{A}^{(h)}_{s}\| \leq \overline{C}^{(h)}
(\overline{\lambda}^{(h)})^{[s-1]},$ nous avons l'estimation
suivante (où $a^{(h)}_{0}=\| \overline{A}^{(h)}_{0} \|$):
\begin{eqnarray*}
\| L^{(h)}_{s} \| &\leq&
\overline{C}^{(h)}\frac{\Gamma(\overline{\lambda}^{(h)}+s-1)}{\Gamma(\overline{\lambda}^{(h)})}+2a^{(h)}_{0}\overline{C}^{(h)}
\frac{\Gamma(a^{(h)}_{0}+s)}{\Gamma(a^{(h)}_{0})}\\&& \ \ \
+\overline{C}^{(h)}\sum_{(j,k,l) \in J_{s}
}c_{j,l}^{(k)}\frac{\Gamma(a^{(h)}_{0}+j)}{\Gamma(a^{(h)}_{0})}\frac{\Gamma(\overline{\lambda}^{(h)}+l-1)}{\Gamma(\overline{\lambda}^{(h)})}
\end{eqnarray*}
qui s'écrit (voir \cite{duvseriesfact}):
\begin{eqnarray*}
&&\overline{C}^{(h)}
\frac{1-\overline{\lambda}^{(h)}}{a^{(h)}_{0}+1-\overline{\lambda}^{(h)}}
\frac{\Gamma(\overline{\lambda}^{(h)}+s-1)}{\Gamma(\overline{\lambda}^{(h)})} \\
 &+&\left(
2a^{(h)}_{0}+\frac{\overline{C}^{(h)}}{a^{(h)}_{0}+1-\overline{\lambda}^{(h)}}\right)
\frac{\Gamma(a^{(h)}_{0}+s)}{\Gamma(a^{(h)}_{0})}
\end{eqnarray*}
et qui est majoré par
$3\overline{C}^{(h)}(\overline{\lambda}^{(h)})^{[s-1]}$ pour $h$
suffisamment petit, uniformément en $s.$ D'autre part, le spectre de l'opérateur $L^{(h)}_{0}$, qui est égal à $\{\mu - \nu \ | \ \mu, \nu \in Sp(\overline{A}^{(h)}_{0})\}$, ne contient aucun entier strictement négatif puisque $\overline{A}^{(h)}_{0}$ est non résonnante. De plus, $I$ appartient au noyau de $L^{(h)}_{0}$.  Nous sommes donc en
mesure d'appliquer la proposition \ref{preparatif}; cela termine la démonstration du théorème.

\subsection{Bilan: la solution canonique en $+\infty$.} Nous sommes à présent en mesure
d'associer une solution canonique à une famille de systèmes
fuchsiens $(C,\lambda)$ en $+\infty,$ non résonnants.

\begin{theo}\label{sol can dans le cas non resonnant}
Soit $\delta_{-h}Y=A^{(h)}Y$, $h\in\io$ une famille de systèmes
fuchsiens $(C,\lambda)$ en $+\infty,$ non résonnants. On note
$A^{(h)}_{0}$ la valeur en $+\infty$ de $A^{(h)}.$ Nous appelons
solutions canoniques, en $+\infty$, de cette famille de systèmes, la famille de
(systèmes fondamentaux de) solutions :
\begin{equation*}
e_{A^{(h)}}^{(h)+}:=F^{(h)+}_{A^{(h)}}e_{A_{0}^{(h)}}^{(h)+}
\end{equation*}
où $e_{A_{0}}^{(h)+}$ a été définie en \ref{resolution des systeme
a matrice constante} et $F^{(h)+}_{A^{(h)}}$ est l'unique transformation de jauge de 
$Gl_{n}(\seriesfact^{(h)}),$ tangente à $I$ en $+\infty,$ telle que $(A^{(h)})^{F^{(h)+}_{A^{(h)}}}=A^{(h)}_{0}$. 
Si on suppose de plus que :
\begin{equation*}
\sup_{s\in \integer^{*}, \ h\in\io}|||K_{s,A_{0}^{(h)}}^{-1}|||<\infty
\end{equation*}
alors la famille des $F^{(h)+}_{A^{(h)}}$, $h\in\io$ est de type
$(C',\lambda')$ pour certains $C',\lambda'$.
\end{theo}

\subsection{Convergence des séries de
factorielles vers les séries entières à l'infini.}

Il est naturel de s'intéresser à la convergence des séries de
factorielles vers les séries entières à l'infini puisque
$x^{-[s]_{h}} \xrightarrow[h \longrightarrow 0]{} x^{-s}$.

\begin{prop} \label{conf de series de fact vers serie entiere en l'infini}
Soit, pour tout $h\in\io,$ $ Y^{(h)}(x)=\sum_{s=0}^{+\infty} Y^{(h)}_{s}x^{-[s]_{h}}$ un élément de $M_n(\seriesfactform^{(h)})$.
Soit également, $ \widetilde{Y}(x)=\sum_{s=0}^{+\infty}
\widetilde{Y}_{s}x^{-s}$ un élément de $M_n(\complex[[x^{-1}]])$. Si
$Y^{(h)}_{s} \xrightarrow[h \longrightarrow 0^+]{}
\widetilde{Y}_{s}$ et si la famille des $Y^{(h)}$, $h\in \io$ est de type
$(C,\lambda)$ alors, pour $h>0$ assez petit, la série de $h$-factorielles formelle
définissant $Y^{(h)}$ converge pour $Re(x)>\lambda+1;$ la série
formelle définissant $\widetilde{Y}$ converge absolument pour
$|x|>\lambda+1$ et $Y^{(h)}(x) \xrightarrow[h \longrightarrow
0^{+}]{} \widetilde{Y}(x)$ pour $Re(x)>\lambda+1.$ La convergence
est uniforme sur tout domaine de la forme $Re(x)>\lambda
+1+\epsilon$ (pour tout $\epsilon > 0$).
\end{prop}

\begin{proof} Commençons par la convergence simple.
Soit $Re(x)>\lambda + 1.$ On suppose $h>0$ assez petit pour que les estimations de type $(C,\lambda)$ soient valables. Nous avons :

\begin{equation*}
 |Y^{(h)}_{s}x^{-[s]_{h}}|\leq C
 \left|\frac{\lambda^{[s-1]_{h}}}{x^{[s]_{h}}}\right| \leq C
\frac{\lambda^{[s-1]_{h}}}{Re(x)^{[s]_{h}}}.
\end{equation*}
De l'inégalité $\frac{1}{Re(x)+(s-1)h}\leq \frac{1}{Re(x)}$ et de la croissance de la fonction $h \longmapsto \frac{\lambda+hk}{Re(x)+hk}$, on déduit :
\begin{align*}
 |Y^{(h)}_{s}x^{-[s]_{h}}|&\leq C
 \left[\prod_{k=0}^{s-2}\left[\frac{\lambda+k}{Re(x)+k}\right]\right]
\frac{1}{Re(x)}.
\end{align*}
La formule de Stirling permet de conclure à la convergence de la
série de terme général
$\prod_{k=0}^{s-2}\left[\frac{\lambda+k}{Re(x)+k}\right].$ Le
théorème de convergence dominée permet de terminer la preuve de la première partie de la
proposition.

Passons à la convergence uniforme.
Soient $\nu>\lambda + 1$ et $Re(x)>\nu + \delta$ ($\delta > 0$).
Effectuons la transformation d'Abel suivante :

\begin{equation*}
  b_{s}^{(h)}=a_{s}^{(h)}\nu^{-[s]_{h}}, \ \
  c_{s}^{(h)}=\nu^{[s]_{h}}x^{-[s]_{h}}, \ \
  B_{s}^{(h)}=\sum_{s=m}^{n}b_{k}^{(h)},
\end{equation*}

\begin{eqnarray*}
  \sum_{s=m}^{n}a_{s}^{(h)}\nu^{-[s]_{h}}&=&\sum_{s=m}^{n-1}B_{s}^{(h)}(c_{s}^{(h)}-c_{s+1}^{(h)}) +
  c_{n}^{(h)}.
\end{eqnarray*}
Nous avons:
\begin{eqnarray*}
 c_{s}^{(h)}-c_{s+1}^{(h)}&=&c_{s}^{(h)}\frac{x-\nu}{x+h(s-1)}\\
 &=&\underbrace{\frac{\nu (x-\nu)}{x}}_{(1)}  \underbrace{\left[\prod_{k=1}^{s-1} \frac{\nu
 +hk}{x+hk} \right]\frac{1}{x+hs}}_{(2)}.
\end{eqnarray*}
D'abord : $(1) \leq 2 |\nu|$,
ensuite : $(2)\leq \frac{ d_{s}^{(h)}}{\nu + \delta + sh}$
avec:
\begin{eqnarray*}
 d_{s}^{(h)}=\prod_{k=1}^{s-1} \frac{\nu
 +hk}{\nu+\delta+hk} \geq 0
\end{eqnarray*}
or
$d_{s}^{(h)}-d_{s+1}^{(h)}=d_{s}^{(h)}\frac{\delta}{\nu+\delta+hs}$,
donc :
\begin{eqnarray*}
 |c_{s}^{(h)}-c_{s+1}^{(h)}| \leq 2 |\nu| \delta^{-1}
 (d_{s}^{(h)}-d_{s+1}^{(h)}).
\end{eqnarray*}
Il apparaît d'autre part clairement que, $\epsilon$ étant donné,
si $N$ est suffisamment grand, alors pour tout $n>m\geq N,$
$|B_{s}^{(h)}|\leq \frac{\epsilon}{2 |\nu| \delta^{-1}}.$ Il
s'en suit que, pour $n>m\geq N$ :
\begin{eqnarray*}
  \sum_{s=m}^{n}a_{s}^{(h)}x^{-[s]_{h}}&\leq&\epsilon d_{s}^{(h)}\leq \epsilon, \\
\end{eqnarray*}
ce qui achève la preuve.
\end{proof}

Le résultat suivant (dont l'utilité apparaîtra clairement au
corollaire \ref{ca c est un super resultat} de la section \ref{etude globale}) fournit une
réciproque partielle à la proposition précédente.

\begin{prop} \label{theo conf unif}
Soit $\widetilde{A}$ un élément de $M_n(\complex(x))$ qu'on
suppose holomorphe en $\infty$. Soit $A^{(h)}$, $h\in \io$ une famille
d'éléments de $M_{n}(\complex(x))$ qui convergent uniformément
vers $\widetilde{A}$ sur un voisinage de $\infty$, lorsque $h$ tend vers $0^+$ (en particulier
les $A^{(h)}$ sont holomorphes en $\infty$ pour $h>0$ assez
petit). Notons :
$$\widetilde{A}(x)=\sum_{s=0}^{+\infty} \widetilde{A}_{s}x^{-s}$$
le développement en série entière au voisinage de l'infini de
$\widetilde{A}$ et : $$A^{(h)}(x)=\sum_{s=0}^{+\infty}
A^{(h)}_{s}x^{-[s]_h}$$ le développement en série de factorielles
de $A^{(h)}.$ Alors la famille de séries de $h$-factorielles $A^{(h)}$,
$h\in \io$ est de type $(C,\lambda)$ et on a, pour tout $s\in \integer$
:
$$A^{(h)}_{s}
\xrightarrow[h \longrightarrow 0^{+}]{} \widetilde{A}_{s}.$$
\end{prop}

\begin{proof}
La seconde assertion pourrait être justifiée par un raffinement de
la méthode utilisée pour prouver la première (\textit{i.e.}
effectuer les calculs exacts dans le lemme \ref{a le beau lemme})
mais nous préférons utiliser un résultat de \cite{watson}  pp. 142-143. On pose
:

\begin{equation*}
  \varphi^{(h)}(x)=\sum_{i=0}^{+\infty}\frac{\underline{A}_{i}^{(h)}}{h^{i}}\frac{(-\ln(1-x))^{i}}{i!}
\text{ et }
  \psi^{(h)}(x)=\varphi^{(h)}(hx).
\end{equation*}
Alors, selon \cite{watson},
$A_{i}^{(h)}=h^{i}\derivnieme{\varphi^{(h)}}{x}{i}(0)=\derivnieme{\psi^{(h)}}{x}{i}(0).$
On constate que, sur un certain voisinage de $0,$ $\psi^{(h)}(x)$
(qui est holomorphe sur un voisinage fixe de $0$) converge
uniformément vers
$\widetilde{\psi}(x)=\sum_{i=0}^{+\infty}\underline{A}_{i}^{(h)}\frac{x^{i}}{i!}.$
Il en résulte que $\lim_{h \longrightarrow 0^{+}} A_{i}^{(h)} =
\widetilde{A}_{i}.$ Ce qui prouve la seconde assertion.

Prouvons maintenant la première.

\begin{lem}\label{a le beau lemme}
Les coefficients du développement en série de
\textit{h}-factorielles de $\frac{1}{x^{n}}=\sum_{s=n}^{+\infty}
\psi_{n,s}^{(h)} x^{-[s]_{h}}$ sont tous positifs.
\end{lem}
\begin{proof}
On commence par le cas $h=1.$ On procède par récur\-rence sur $n.$
Le résultat est trivial pour $n=0$ (ainsi que pour $n=1$).
Supposons alors que l'hypothèse est vérifiée au rang $n$ :
\begin{equation*}
\frac{1}{x^{n}}=\sum_{s=n}^{+\infty} \psi_{n,s}^{(1)} x^{-[s]},\
\psi_{n,s}^{(1)}\geq 0.
\end{equation*}
Dérivant les deux côtés de cette égalité, nous obtenons, grâce à
la formule de dérivation des séries de factorielles donnée par Norlund dans \cite{Norlund
series d interp} pp. 220-222:
\begin{equation*}
-\frac{n}{x^{n+1}}=-\sum_{s=n+1}^{+\infty}
\left(\frac{\psi_{n,s-1}^{(1)}}{1(s-2)!}+\frac{\psi_{n,s-2}^{(1)}}{2(s-3)!}+...+\frac{\psi_{n,1}^{(1)}}{s-1}\right)(s-1)!
x^{-[s]},
\end{equation*}
donc :
\begin{equation*}
\psi_{n+1,s}^{(1)}=\left(\frac{\psi_{n,s-1}^{(1)}}{1(s-2)!}+\frac{\psi_{n,s-2}^{(1)}}{2(s-3)!}+...+\frac{\psi_{n,1}^{(1)}}{s-1}\right)(s-1)!/n\geq
0
\end{equation*}
et la récurrence est achevée.

Pour le cas $h$ quelconque, il suffit de substituer $x/h$ à $x$
pour obtenir le résultat.
\end{proof}

\begin{lem}
Nous avons :
$\sum_{k=0}^{+\infty} \psi_{k,s}^{(h)}
\alpha^{k-1}=\alpha^{[s-1]_{h}}.$
\end{lem}
\begin{proof}
Nous avons les développements suivants :
\begin{equation*}
\frac{1}{x-\alpha}=\sum_{s=0}^{+\infty}\alpha^{[s-1]_{h}}x^{-[s]_{h}}=\sum_{s=0}^{+\infty}\alpha^{s}x^{-(s+1)}.
\end{equation*}
Par identification, nous en déduisons :
$\sum_{k=0}^{+\infty} \psi_{k,s}^{(h)}
\alpha^{k-1}=\alpha^{[s-1]_{h}},$
ce qui est le résultat cherché.
\end{proof}

Par hypothèse, il existe un voisinage de l'infini $D=\{|x|\geq
r\}$ sur lequel les $A^{(h)}$ (n'ont pas de pôles et) convergent
uniformément vers $\widetilde{A}.$ Il résulte alors des
estimations de Cauchy, que pour $h>0$ assez petit, il existe une constante $C>0$ telle que, pour tout $i\in\integer$:
\begin{equation*}
\|A_{i}^{(h)}\|\leq Cr^{i}.
\end{equation*}
Avec les notations ci-dessus nous avons:
\begin{equation*}
A_{i}^{(h)}=\sum_{k=0}^{+\infty} \psi_{k,i}^{(h)}
\underline{A}_{k}^{(h)}
\end{equation*}
et, compte tenu de ce que $\psi_{k,i}^{(h)}\geq 0$ pout tout $k,\
i,$ on en déduit que, pour $h>0$ assez petit et pour tout $i\in \integer$ :
\begin{equation*}
\|A_{i}^{(h)}\|\leq C\sum_{k=0}^{+\infty} \psi_{k,i}^{(h)}
r^{i}=(Cr) r^{[i-1]_{h}}.
\end{equation*}
Autrement dit, la famille des $A^{(h)}$ est de type $(Cr,r).$ Ceci vient clore la preuve du théorème.
\end{proof}

\subsection{Confluence de différence vers différentielle.}  Pour appliquer la méthode dite de Frobenius à un système différentiel, on emploiera
la détermination principale du logarithme. Par souci de
simplicité, on opte ici pour les caractères et les logarithmes de
l'exemple 1. Il n'y a aucune difficulté à adapter les énoncés pour
les choix de l'exemple 2.

\subsubsection{Etude locale au voisinage de $+\infty.$}

Nous introduisons l'hypothèse suivante :

\begin{description}
  \item[(H)] (i) Soit, pour tout $h\in\io$, $A^{(h)} \in M_n(\seriesfact^{(h)})$. On suppose que la famille $A^{(h)}$, $h\in\io$ est de type $(C,\lambda)$, et que, pour tout $s\in \integer,$ il existe $\widetilde{A}_{s} \in M_n(\complex)$ avec
$A^{(h)}_{s} \xrightarrow[h \longrightarrow 0^{+}]{}
\widetilde{A}_{s}$ (où $A^{(h)}(x)=\sum_{s=0}^{+\infty}
A^{(h)}_{s}x^{-[s]_h}$).\\

(ii) On suppose qu'une réduction de Jordan de $\widetilde{A}_{0}$
se
déploie en des réductions de Jordan des $A^{(h)}_{0}$, $h\in \io$ (\textit{cf.} définition \ref{def deploiement}). \\

(iii) On suppose que $\widetilde{A}_0$ est non résonnante.
\end{description}

On suppose que (H) est satisfaite.

D'après la proposition \ref{conf de series de fact vers serie entiere en l'infini} la série formelle
$\widetilde{A}(x)=\sum_{s=0}^{+\infty} \widetilde{A}_{s}x^{-s}$
est convergente au voisinage de l'infini. Ainsi le système différentiel :
 $$\delta
\widetilde{Y}=\widetilde{A}\widetilde{Y}$$ est fuchsien en $\infty$
et non résonnant. On note $e_{\widetilde{A}}$ sa solution canonique,
en $\infty$, obtenue par la méthode de Frobenius. 

D'autre part, remarquons que les systèmes
(aux différences) définis par les $A^{(h)}$ sont non résonnants pour
$h>0$ assez petit : cela résulte immédiatement de l'hypothèse de
déploiement des réductions de Jordan et de la non
résonnance de $\widetilde{A}_0.$

\begin{theo}
 Si (H) est vérifiée alors $e_{A^{(h)}}^{(h)+}$ tend, quand
$h$ tend vers $0^{+},$ sur un certain demi-plan droit, vers $e_{\widetilde{A}}$. La
convergence est uniforme sur tout demi-plan droit $Re(x)>C,$ pour
$C$ assez grand.
\end{theo}

\begin{proof}
 Notons que l'hypothèse
$\sup_{s\in \integer^{*}, \ h\in ]0,h_0'[}|||K_{s,A_{0}^{(h)}}^{-1}|||<\infty$
du théorème \ref{sol can dans le cas non resonnant} est ici
vérifiée (pour $h'_0\in \io$ assez petit). En effet, on a :

\begin{equation*}
K_{s,A_{0}^{(h)}}(U) =
\underbrace{\left[\widetilde{A}_{0}U-U\widetilde{A}_{0}-sU\right]}_{=:B_{s}^{-1}(U)}+\underbrace{\left[(A^{(h)}_{0}-\widetilde{A}_{0})U-U(A^{(h)}_{0}-\widetilde{A}_{0})\right]}_{=:C^{(h)}(U)},
\end{equation*}
($B_{s}$ et $C^{(h)}$ sont des applications linéaires) ainsi:

\begin{equation*}
 |||K_{s,A_{0}^{(h)}}^{-1}|||\leq
a\sum_{k=0}^{+\infty} (bc^{(h)})^{k}
\end{equation*}
où $a=\sup_{s\in \integer^{*}, \ h>0}|||B_{s}^{-1}|||<\infty$, $b=\sup_{s\in \integer^{*}, \ h>0}|||B_{s}|||<\infty$ et
$c^{(h)}=||| C^{(h)} |||$ mais $c^{(h)}$ tend vers $0$ avec $h$
(uniformément en $s\in
\integer^{*}$) et on a justifié l'assertion.\\

On a donc la famille de solutions canoniques
$e_{A^{(h)}}^{(h)+}:=F^{(h)+}_{A^{(h)}}e_{A_{0}^{(h)}}^{(h)+};$ on
abrège $F^{(h)+}_{A^{(h)}}$ en $F^{(h)}.$ On traite séparément la
partie ``log-car'' et la partie ``transformation de jauge''.\\

La partie ``log-car" a déjà été étudiée :
c'était l'objet de la proposition \ref{conflu dans le cas des coeff constants}.\\

Passons à la partie ``transformation de jauge''. On note
$\widetilde{F}$ la transformation de jauge intervenant dans la
résolution de $\delta \widetilde{Y}=\widetilde{A}\widetilde{Y}.$
D'après le théorème \ref{sol can dans le cas non resonnant}, la famille $F^{(h)}$, $h\in ]0,h'_0[$ est de type $(C,\lambda);$ ainsi, si on
arrive à prouver que $F^{(h)}_{s} \longrightarrow
\widetilde{F}_{s},$ on pourra conclure grace à la proposition
\ref{conf de series de fact vers serie entiere en l'infini}.

La transformation de jauge $F^{(h)}$ est caractérisée (en reprenant les notations introduites en \ref{qsd}) par :
$\overline{F}^{(h)}$ est une fonction développable en série de
factorielles tangente à l'identité en l'infini telle que
$\overline{A}^{(h)}(x)\overline{F}^{(h)}(x)-\delta_{-1}
\overline{F}^{(h)}(x)=\overline{F}^{(h)}(x-1)\overline{A}^{(h)}_{0}.$

Or:
\begin{equation*}
\delta_{-1}\overline{F}^{(h)}(x)=\sum_{s=0}^{+\infty}
-s\overline{F}^{(h)}_{s}x^{-[s]}
\end{equation*}
et (formule de translation; voir par exemple \cite{duvseriesfact}):
\begin{equation*}
\overline{F}^{(h)}(x-1)=\overline{F}^{(h)}_{0}+\sum_{s=1}^{+\infty}\left[
\overline{F}^{(h)}_{s}+(s-1)!\sum_{k=1}^{s-1}
\frac{\overline{F}^{(h)}_{k}}{(k-1)!}\right]x^{-[s]}
\end{equation*}
puis :
\begin{equation*}
\overline{A}^{(h)}(x)\overline{F}^{(h)}(x)=\sum_{s=0}^{+\infty}C_{s}x^{-[s]}
\end{equation*}
avec:
\begin{equation*}
C_{0}=\overline{A}^{(h)}_{0}\overline{F}^{(h)}_{0}, \
C_{s}=\overline{A}^{(h)}_{0}\overline{F}^{(h)}_{s}+\overline{A}^{(h)}_{s}\overline{F}^{(h)}_{0}+\sum_{(j,k,l)
\in J_{s}
}c_{j,l}^{(k)}\overline{A}^{(h)}_{j}\overline{F}^{(h)}_{l}
\end{equation*}
et:
\begin{equation*}
J_{s} = \{(j,k,l) \backslash j,l \geq 1,\ k \geq 0,\ j+k+l=s\}
\end{equation*}
\begin{equation*}
c_{j,l}^{(k)}=\frac{(j+k-1)!(l+k-1)!}{k!(j-1)!(l-1)!}.
\end{equation*}
On en déduit que $F^{(h)}$ vérifie :

  \begin{equation}\label{super}\begin{cases}
    \overline{F}^{(h)}_{0}=I \\
    \phi_{\overline{A}^{(h)}_{0},s}(\overline{F}^{(h)}_{s})=-\overline{A}^{(h)}_{s}+(s-1)!\left[\sum_{k=1}^{s-1}
    \frac{\overline{F}^{(h)}_{k}}{(k-1)!}\right]\overline{A}^{(h)}_{0}  \\
    \ \ \  \ \ \ \ \ \ \ \ \ \ \ \ \ \  \ \ \ \ \ \ \ \ \ \ \ \ \  \ \ \ \ \ \ \ \ \ \ \ -\sum_{(j,k,l) \in
J_{s} }c_{j,l}^{(k)}\overline{A}^{(h)}_{j}\overline{F}^{(h)}_{l}
  \end{cases}\end{equation}
où :
\begin{eqnarray*}
\phi_{\overline{A}^{(h)}_{0},s} : M_{n}(\complex) & \longrightarrow &
M_{n}(\complex)\\
U & \longmapsto & (\overline{A}^{(h)}_{0}+sI)U-U\overline{A}^{(h)}_{0}.
\end{eqnarray*}

Il est bien connu que si $\overline{A}^{(h)}_{0}$ est non résonnante alors $\phi_{\overline{A}^{(h)}_{0},s} :
M_{n}(\complex)  \longrightarrow  M_{n}(\complex)$ est un
isomorphisme pour tout $s \neq 0.$ La non résonnance permet donc
de résoudre de manière unique les
 relations de récurrence (\ref{super}) : $F^{(h)}$ est caractérisée par :

 \begin{equation*} \begin{cases}
    \overline{F}^{(h)}_{0}=I \\
    \overline{F}^{(h)}_{s}=\phi_{\overline{A}^{(h)}_{0},s}^{-1}(-\overline{A}^{(h)}_{s}+(s-1)!\left[\sum_{k=1}^{s-1}
    \frac{\overline{F}^{(h)}_{k}}{(k-1)!}\right]\overline{A}^{(h)}_{0}\\
   \ \ \ \ \ \ \ \ \ \ \ \ \ \ \ \ \ \ \ \ \ \ \ \ \ \ \ \ \ \ \ \ \ \ \ \ \ -\sum_{(j,k,l) \in
J_{s}
}c_{j,l}^{(k)}\overline{A}^{(h)}_{j}\overline{F}^{(h)}_{l}),\ s
\in \integer^{*}.
  \end{cases}\end{equation*}
Maintenant, par récurrence, on voit que $F^{(h)}_{s}$ admet une
limite pour $h \longrightarrow 0^+$ qu'on note
$\widetilde{F}'_{s}.$ Il est clair que les $\widetilde{F}'_{s}$
satisfont les relations de récurrence des coefficients de la
transformation de jauge du système différentiel limite (non
résonnant), qui intervient dans la méthode de Frobenius. Ainsi
$\widetilde{F}'_{s}=\widetilde{F}_{s}.$
\end{proof}

\subsubsection{Etude globale.} \label{etude globale} On se place à présent dans le cas
algébrique : \textbf{les matrices qui définissent les systèmes aux
différences et différentiels sont supposées rationnelles}. On va
``tirer'' le demi-plan droit introduit pour la transformation de
jauge précédente vers la gauche. Pour ce faire, on ajoute
l'hypothèse suivante :\\

\begin{description}
  \item[(H')] (i) Soient $\widetilde{A}$ un élément de
  $M_{n}(\complex(x))$ et $A^{(h)}$, $h\in \io$ une famille d'éléments de $M_{n}(\complex(x))$ qui vérifient
  (H).

  (ii) On suppose que $\widetilde{A}$ est la limite
uniforme, lorsque $h$ tend vers $0^+$, des $A^{(h)}$ sur tout
compact de $\complex \setminus
\{\text{pôles de $\widetilde{A}$}\}$.\\
\end{description}

Dans le point (ii) de (H'), on sous-entend que, étant donné
$x_{0}$ un nombre complexe qui n'est pas un pôle de
$\widetilde{A},$ il existe un voisinage de $x_0$ sur lequel les
$A^{(h)}$ n'ont pas de pôle, pour $h>0$ assez petit.\\

\textbf{Notations. } On note $\widetilde{\Omega}^+$ le plan
complexe privé des demi-droites horizontales issues de chacun des
pôles de
$\widetilde{A}$ et de $0$ et dirigées vers $-\infty.$\\

On suppose (H') vérifiée. Rappelons qu'alors le système
différentiel $\delta \widetilde{Y} = \widetilde{A} \widetilde{Y}$
est fuchsien en $\infty$ et non résonnant. On note
$e_{\widetilde{A}}$ le prolongement à $\widetilde{\Omega}^+$ de sa
solution
canonique, en $\infty$, obtenue par la méthode de Frobenius.

Rappelons également que les systèmes
(aux différences) définis par les $A^{(h)}$ sont non résonnants pour
$h>0$ assez petit.\\

-\textbf{Préliminaires}

Nous utiliserons les deux résultats suivants. Le premier est issu de \cite{sujet
agreg} alors que le second provient de \cite{sauloy  systemes aux q diff
singu reg }.

\begin{theo}
Soient $\mathcal{A}$ une fonction continue sur un segment réel
$]a,b[$ à valeurs matricielles complexes et $t_{0}<t_{1}$ dans
$]a,b[.$ Soit $p$ un entier non nul et
$s_{0}^{(p)}<\cdots<s_{p}^{(p)}$ une subdivision de
$]t_{0},t_{1}[$ de pas (maximum de la différence de deux éléments
successifs de la subdivision) $\epsilon_{p}.$ On suppose que
$\epsilon_{p}$ tend vers $0$ quand $p$ tend vers $+\infty$. Alors la
résolvante du système $\mathcal{X}'=\mathcal{A}\mathcal{X}$ est
donnée par :

\begin{equation*}
 R(t_{1},t_{0})=\lim_{p\longrightarrow \infty}
 (I+\mathcal{A}(s^{(p)}_{p-1})\Delta s^{(p)}_{p-1})\cdots (I+\mathcal{A}(s^{(p)}_{0})\Delta s^{(p)}_{0})
\end{equation*}
avec $\Delta s^{(p)}_{i}=s^{(p)}_{i+1}-s^{(p)}_{i}.$
\end{theo}

\begin{prop}
Soient $(A_{p,i})$ et $(B_{p,i})$ deux familles de matrices
indexées par les couples d'entiers $(p,i)$ tels que $1\leq i \leq p.$ On
considère les suites de termes généraux
$R_{p}=(I+A_{p,1})\cdots(I+A_{p,p})$ et
$S_{p}=(I+B_{p,1})\cdots(I+B_{p,p}).$ On suppose que, quand $p
\longrightarrow +\infty$ :
\begin{itemize}
  \item la suite $(\Sigma_{1\leq i \leq p}\|A_{p,i}\|)$ est
  bornée;
  \item la suite $(\Sigma_{1\leq i \leq p}\|A_{p,i}-B_{p,i}\|)$
  tend vers $0.$
\end{itemize}
Alors, $\lim_{p\longrightarrow \infty} \|R_{p}-S_{p}\|=0.$
\end{prop}

-\textbf{Etude}

Notre résultat ``global'' est le suivant.

\begin{theo}
Supposons (H'). Alors $e_{A^{(h)}}^{(h)+}$ tend, sur
$\widetilde{\Omega}^+$, quand $h$ tend vers $0^{+}$, vers
$e_{\widetilde{A}}.$
\end{theo}

\begin{proof}
Afin d'alléger la présentation, on pose :
$$Y^{(h)}=e_{A^{(h)}}^{(h)+} \text{ et } \widetilde{Y}=e_{\widetilde{A}}.$$
Nous introduisons les notations suivantes où $t,\ \eta \in \real$ :
\begin{equation*}
\check{Y}^{(h)}_{\eta}(t)=Y^{(h)}(-t-i \eta)
\text{ et }
\check{\widetilde{Y}}_{\eta}(t)=\widetilde{Y}(-t-i \eta).
\end{equation*}
 On a:
 \begin{equation*}
\check{Y}^{(h)}_{\eta}(t+h)=(I+\frac{hA^{(h)}(-t-i \eta )}{t+h+i
\eta})\check{Y}^{(h)}_{\eta}(t).
\end{equation*}
Nous savons (étude locale) que, pour $a\ll 0$:
 \begin{equation*}
\check{Y}^{(h)}_{\eta}(t)\xrightarrow[h\longrightarrow 0^{+}
]{}\check{\widetilde{Y}}_{\eta}(t),\ t<a;
\end{equation*}
quitte à choisir $a$ encore plus petit, on peut supposer que les
pôles de $\frac{\widetilde{A}(-x)}{x}$ sont dans le demi-plan
$\{Re(x)>a\}.$ Soient $t_{1}>a>t_{0}$ tels que
$\frac{\widetilde{A}(-t-i\eta)}{t+i\eta}$ n'ait pas de pôle sur
$[t_{0},t_{1}].$ On considère la subdivision suivante de $]t_{0},t_{1}]$ :
 $$t_0 <\nu h< \nu h +h  < \nu +2h <...< \nu h +Nh=t_1$$
où $N=E(\frac{t_{1}-t_{0}}{h})$ et $\nu =
\frac{t_{1}-t_{0}}{h} - N.$ La résolvante $R(t_{0},t_{1})$ du
système :
\begin{equation*}
\frac{d}{dt}
\check{\widetilde{Y}}_{\eta}(t)=\frac{\widetilde{A}(-t-i\eta)}{t+i\eta}\check{\widetilde{Y}}_{\eta}(t)
\end{equation*}
est donnée, compte tenu des préliminaires, par :

\begin{eqnarray*}
R(t_{0},t_{1}) &=&\lim_{h\longrightarrow 0}
 (I+h\frac{\widetilde{A}(h-t_{1}-i\eta)}{t_{1}-h+i\eta})\cdots\\ &&\ \ \ \ \cdots (I+h\frac{\widetilde{A}(Nh-t_{1}-i\eta)}{t_{1}-Nh+i\eta})
(I+\nu h\frac{\widetilde{A}(-t_{0}-i\eta)}{t_{0}+i\eta}).
\end{eqnarray*}

Remarquons que $\frac{A^{(h)}(-t-i \eta )}{t+h+i \eta}$ converge
uniformément vers $\frac{\widetilde{A}(-t-i\eta)}{t+i\eta},$
lorsque $h \longrightarrow 0^{+},$ sur $[t_{0},t_{1}].$ Ainsi,
d'après les préliminaires, la résol\-van\-te s'écrit aussi:

\begin{eqnarray*}
R(t_{0},t_{1}) &=&\lim_{h\longrightarrow 0}
 (I+h\frac{A^{(h)}(h-t_{1}-i\eta)}{t_{1}-h+i\eta})\cdots\\ &&\  \cdots (I+h\frac{A^{(h)}(Nh-t_{1}-i\eta)}{t_{1}-Nh+i\eta})
(I+\nu h\frac{A^{(h)}(-t_{0}-i\eta)}{t_{0}+i\eta}).
\end{eqnarray*}

Par suite :
$\check{Y}^{(h)}_{\eta}(t_{1})\xrightarrow[h\longrightarrow 0^{+}
]{}\check{\widetilde{Y}}_{\eta}(t_{1}).$

\end{proof}

Ce théorème admet un corollaire qui
n'impose pas \textit{a priori} de vérifier le type $(C,\lambda)$ qui est
automatique si on se donne une hypothèse supplémentaire de
convergence uniforme au voisinage de $\infty.$ On introduit à
nouveau une hypothèse.

\begin{description}
  \item[(H'')] (i) Soit $\delta \widetilde{Y}=\widetilde{A}\widetilde{Y}$ un
  système différentiel
algébrique fuchsien en $\infty$ et non résonnant.

(ii) Soient $A^{(h)}$, $h\in\io$ une famille d'éléments de
$M_{n}(\complex(x))$ qui convergent uniformément, lorsque $h$ tend
vers $0^+$, vers $\widetilde{A}$, sur tout compact de
$\mathbb{P}^{1}\complex \backslash \{\text{pôles de
$\widetilde{A}$}\}$.

(iii) On suppose enfin qu'une réduction de Jordan de
$\widetilde{A}_{0}=\widetilde{A}(\infty)$ se déploie en des
réductions de Jordan des $A^{(h)}_{0}$ (\textit{cf.} définition
\ref{def deploiement}).
\end{description}

Dans le point (ii) de (H''), on sous-entend que, étant donné
$x_{0} \in \mathbb{P}^{1}\complex$ qui n'est pas un pôle de
$\widetilde{A}$, il existe un voisinage de $x_0$ sur lequel les
$A^{(h)}$ n'ont pas de pôle, pour $h>0$ assez petit; en
particulier $\infty$ est un point régulier des $A^{(h)}$ pour
$h>0$ assez petit.

On note à nouveau $e_{\widetilde{A}}$ le prolongement à
$\widetilde{\Omega}^+$ de la solution canonique, en $\infty$, du
système algébrique fuchsien en $\infty$ et non résonnant $\delta
\widetilde{Y} = \widetilde{A} \widetilde{Y}$, obtenue par la
méthode de Frobenius.

\begin{coro} \label{ca c est un super resultat}
 Si (H'') est vérifiée alors, pour $h'_0\in ]0,h_0[$ assez petit, la
famille de systèmes $\delta_{-h}Y=A^{(h)}Y$, $h \in ]0,h_0'[$ est fuchsienne
$(C,\lambda)$ en $+\infty,$ non résonnante et $e_{A^{(h)}}^{(h)+}$
tend, sur $\widetilde{\Omega}^+$, quand $h$ tend vers $0^{+},$
vers $e_{\widetilde{A}}.$
\end{coro}

\begin{proof}
C'est une application du théorème précédent en vertu de la
proposition \ref{theo conf unif}.
\end{proof}

On a des résultats analogues en $-\infty$ par le changement de variables $x \leftarrow -x.$ Lorsqu'elles existent, les solutions canoniques en $-\infty$ seront notées $e_{A^{(h)}}^{(h)-}.$

\subsection{Matrice de connexion et confluence.}

On définit les matrices de connexion de Birkhoff associées à une
famille de systèmes algébriques fuchsiens non résonnants:
\begin{equation*}
\delta_{-h}Y=A^{(h)}Y
\end{equation*}
par:
\begin{equation*}
P^{(h)}_{A^{(h)}}=P^{(h)}=(e_{A^{(h)}}^{(h)+})^{-1}e_{A^{(h)}}^{(h)-}.
\end{equation*}
\\
On introduit une hypothèse :
\begin{description}
  \item[(H''')] (i) Soient $\widetilde{A}$ un élément de
  $M_{n}(\complex(x))$ et $A^{(h)}$, $h\in \io$ une famille d'éléments de $M_{n}(\complex(x))$ qui vérifient
  (H') (ou (H'')) en $+\infty$ et en $-\infty.$

  (ii) On suppose que deux pôles distincts de $\widetilde{A}$ n'ont
  pas la même partie imaginaire et qu'aucun d'entre eux n'est réel.\\
\end{description}

On suppose maintenant (H''').\\

 \textbf{Notations. } On note $\mathcal{P}(\widetilde{A})$
 l'ensemble des
 pôles de $\widetilde{A}$ et on pose $\mathcal{P}(\widetilde{A})\bigcup
 \{0\}=\cup_{j=1}^{r}\{\widetilde{z}_{j}\}$ où les
 $\widetilde{z}_{j}$ sont deux à deux distincts et rangés selon
 l'ordre croissant de leurs parties imaginaires.\\

 On définit $\widetilde{\Omega}=\complex \setminus
\bigcup_{j=1}^{r}\left(\widetilde{z}_{j}+\real\right)$ ainsi que
$\widetilde{\Omega}^{+}=\complex \setminus
\bigcup_{j=1}^{r}\left(\widetilde{z}_{j}+\real^{-}\right)$ et
$\widetilde{\Omega}^{-}=\complex \setminus
\bigcup_{j=1}^{r}\left(\widetilde{z}_{j}+\real^{+}\right).$
On a donc $\widetilde{\Omega}=\widetilde{\Omega}^{+}\bigcap \widetilde{\Omega}^{-}.$\\

On note $e_{\widetilde{A}}^{+}$ (resp. $e_{\widetilde{A}}^{-}$) la
limite de $e_{A^{(h)}}^{(h)+}$ (resp. $e_{A^{(h)}}^{(h)-}$), sur
$\widetilde{\Omega}^{+}$ (resp. $\widetilde{\Omega}^{-}$), lorsque
$h \longrightarrow 0^{+}$ (dont l'existence est assurée par les
résultats de la section précédente). Alors $P^{(h)} \xrightarrow[h
\longrightarrow 0^{+}]{}
\widetilde{P}=(e_{\widetilde{A}}^{+})^{-1}e_{\widetilde{A}}^{-}$
sur $\widetilde{\Omega}.$ Puisque $e_{\widetilde{A}}^{+}$ et
$e_{\widetilde{A}}^{-}$ sont des systèmes fondamentaux de
solutions de la même équation différentielle, nous avons $\delta
\widetilde{P}=0:$ $\widetilde{P}$ est constante sur les
composantes connexes de $\widetilde{\Omega};$ on note
$\widetilde{P}_{j+1}$ la valeur de $\widetilde{P}$ sur la
composante connexe de $\widetilde{\Omega}$
ayant $\widetilde{z}_{j}+\real$ comme droite frontière inférieure
et $\widetilde{P}_{1}$ sa valeur sur la composante restante. On a alors le

\begin{theo}
Si (H''') est vérifiée alors la monodromie du système différentiel
algébrique (fuchsien en $\infty$) $\delta
\widetilde{Y}=\widetilde{A}\widetilde{Y}$ autour de
$\widetilde{z}_{j}$ dans la base $e_{\widetilde{A}}^{+}$ est
$\widetilde{P}_{j}\widetilde{P}_{j+1}^{-1}.$
\end{theo}

\begin{proof}
En effet, un cercle direct $\gamma_{j}$ (de rayon assez petit)
autour de $\widetilde{z}_{j}$ peut être vu comme un chemin
$\gamma_{j}^{+}$ dans $\widetilde{\Omega}^{+}$ suivi d'un chemin
$\gamma_{j}^{-}$ dans $\widetilde{\Omega}^{-}.$ Le prolongement
analytique le long de $\gamma_{j}^{+}$ transforme
$e_{\widetilde{A}}^{+}$ en
$e_{\widetilde{A}}^{-}\widetilde{P}_{j+1}^{-1};$ celui le long de
$\gamma_{j}^{-}$ transforme $e_{\widetilde{A}}^{-}$ en
$e_{\widetilde{A}}^{+}\widetilde{P}_{j}$ et donc, le long de
$\gamma_{j},$ $e_{\widetilde{A}}^{+}$ est transformé en
$e_{\widetilde{A}}^{-}\widetilde{P}_{j}\widetilde{P}_{j+1}^{-1}.$
\end{proof}

Les hypothèses de ce théorème n'excluent pas la
présence de singularités irrégulières à distance finie de
l'origine.

\subsection{Exemples.}

\subsubsection{$A^{(h)}(x)=A\left(\frac{x}{h}\right),$ $A\in M_{n}(\complex(x))$ holomorphe à l'infini et non résonnante.}

Il est aisé de vérifier que la famille des systèmes définis par
$A^{(h)}(x)=A\left(\frac{x}{h}\right)$
rentre dans le cadre de notre étude et donc que la solution
canonique en $+\infty$ tend, quand $h$ tend vers $0^{+},$ vers
$x^{A_{0}}$ sur $\complex \setminus \real^-.$

\subsubsection{$A^{(h)}(x)=A(x),$ $A\in M_{n}(\complex(x))$ holomorphe à l'infini et non résonnante.}
\label{ex de def le plus simple}

Puisque $A^{(h)}(x)=A(x)$ converge uniformément vers $A(x)$ sur
tout compact de $\mathbb{P}^{1}\complex$ privé des pôles de $A,$
cet exemple fait partie du champ d'application de nos résultats.

\subsubsection{Un exemple régulier en dimension 1:
$\delta_{-h}y(x)=\frac{-\mu}{x-h-\lambda}y(x).$}

La famille des systèmes
$\delta_{-h}y(x)=\frac{-\mu}{x-h-\lambda}y(x)$ avec $Im(\lambda)>0$ vérifie les
hypothèses du corollaire \ref{ca c est un super resultat}. Les
solutions canoniques sont données (en $+\infty$ puis en $-\infty$)
par:
\begin{equation*}
\frac{\Gamma(\frac{x}{h})\Gamma(\frac{x-\lambda}{h})}{\Gamma(\frac{x}{h}-\alpha_{1}^{(h)})
\Gamma(\frac{x}{h}-\alpha_{2}^{(h)})}
\text{ et }
\frac{\Gamma(-\frac{x}{h}+1+\alpha_{1}^{(h)})\Gamma(-\frac{x}{h}+1+\alpha_{2}^{(h)})}{\Gamma(-\frac{x}{h}+1)
\Gamma(-\frac{x}{h}+1+\frac{\lambda}{h})}
\end{equation*}
avec:
\begin{equation*}
\alpha_{1}^{(h)}=\frac{\lambda}{h}\left[1-\sqrt{1-4\frac{\mu}{\lambda^{2}}h}\right]
\text{ et }
\alpha_{2}^{(h)}=\frac{\lambda}{h}\left[1+\sqrt{1-4\frac{\mu}{\lambda^{2}}h}\right];
\end{equation*}
on a :
\begin{equation*}
\alpha_{1}^{(h)}\xrightarrow[h\longrightarrow
0^{+}]{}\frac{\mu}{\lambda}
\end{equation*}
et :
\begin{equation*}
\alpha_{2}^{(h)}=\frac{\lambda}{h}\left[1-\frac{\mu}{\lambda^{2}}h+o(h)\right].
\end{equation*}
La formule des compléments permet d'écrire :
\begin{equation*}
P^{(h)}=\frac{\sin(\frac{x}{h}-\alpha_{1}^{(h)})
\sin(\frac{x}{h}-\alpha_{2}^{(h)})}{\sin(\frac{x}{h})\sin(\frac{x-\lambda}{h})}.
\end{equation*}
Par suite:
\begin{equation*}
\widetilde{P}_{1}=1,\
\widetilde{P}_{2}=e^{-2\pi i \frac{\mu}{\lambda}}
\text{ et }
\widetilde{P}_{3}=1.
\end{equation*}
Ainsi, $e^{2 \pi i \frac{\mu}{\lambda}}$ est la monodromie autour
de $0$ et $e^{-2 \pi i \frac{\mu}{\lambda}}$ celle autour de
$\lambda$ (on peut directement vérifier qu'il s'agit du bon
résultat; en effet, la méthode de Frobenius donne, comme solution
canonique du système limite, la fonction $x \longmapsto
\left(\frac{x-\lambda}{x}\right)^{-\frac{\mu}{\lambda}}$ qui a
bien les monodromies ci-dessus).

\appendix

\section{A propos des caractères et logari\-thmes utilisés.}

\subsection{Les ``caractères'' et ``logarithmes'' comme élém\-ents de
$\seriesfact[\log(x)].$}\label{annexe}

D'après \cite{Norlund series d interp}, la fonction $\frac{\Gamma(x)}{\Gamma(x-c)}x^{-c}$
est développable en une série de factorielles convergente et
tangente à $1$ en $+\infty;$ de là nous déduisons, pour $k\geq 1$ :
\begin{equation*}
   x^{-c}\derivpartnieme{}{c}{k}\frac{\Gamma(x)}{\Gamma(x-c)}=\sum_{i=0}^{k-1}\Omega_{i}(x)\log^{i}(x)
   + \Omega_{k}(x)\log^{k}(x)
\end{equation*}
où les $\Omega_{i},$ $i=0,...,n-1$ sont dans $\seriesfact$ et
tangentes à $0$ en $+\infty$ et $\Omega_{n} \in \seriesfact$ et
est tangente à $I$ en $+\infty.$

\subsection{Applications.}
A l'aide de \ref{annexe}, nous pouvons prouver la

\begin{prop}\label{proposition independance des log}
Les logarithmes $l^{(k)}_{c}$, $k \in \integer$ sont linéairement
indé\-pen\-dants sur $\corpsseriesfact.$
\end{prop}

\begin{proof}
Supposons qu'il existe une relation linéaire non triviale $\sum_{k=0}^{n}a_{k}l^{(k)}_{c}=0$
où $a_{k} \in \corpsseriesfact$.
Pour tout entier $k$ tel que $a_{k} \neq 0$, il existe $\alpha_{k}
\in \rinteger$ tel que
$\frac{a_{k}l^{(k)}_{c}}{x^{\alpha_{k}+c}\log^{k}(x)}
\xrightarrow[x \longrightarrow +\infty]{} b_{k} \in \complex^{*}.$
On divise alors la relation linéaire par $x^{k+c}\log^{k}(x)$ où $k$
est le plus grand des
entiers $n$ tels que $\alpha_{n}=\max\{\alpha_{j}\}_{j}.$ On passe
ensuite à la limite pour obtenir une relation de la forme $a=0$
avec $a \in \complex^{*},$ ce qui est absurde.
\end{proof}

\begin{prop}\label{equation quotient de deux caracteres}
L'équation:
\begin{equation} \label{equation du theo de forme
normale} \tau^{-1}y=\frac{x-1-c}{x-1-c'}y
\end{equation} n'admet pas de solution non nulle dans
$\corpsseriesfact(log(x))$ lorsque $d=c-c' \notin \rinteger.$ Si
$d=0,$ les seules solutions dans $\corpsseriesfact(log(x))$ sont
les constantes.
\end{prop}

\begin{proof}
 Soit:
\begin{equation*}
  y=\frac{\sum_{i=0}^{n}log^{i}(x)\Omega_{i}}{\sum_{j=0}^{m}log^{j}(x)\Omega'_{j}}
\end{equation*}
une éventuelle solution de (\ref{equation du theo de forme
normale}) dans $\corpsseriesfact(log(x))$ ($\Omega_{i},\
\Omega'_{j} \in \corpsseriesfact$). Rappelons qu'étant donné $\Omega \in \seriesfact \setminus \{0\}$, il existe $n \in
\rinteger$ tel que $\frac{\Omega}{x^{n}} \in \seriesfact$ et tende
vers une limite non nulle pour $x \rightarrow +\infty$. Ainsi, $y$ peut s'écrire sous la forme :
\begin{equation*}
 y=\frac{\sum_{i=0}^{n}log^{i}(x)x^{\alpha_{i}}\Omega_{i}}{\sum_{j=0}^{m}log^{j}(x)x^{\alpha'_{j}}\Omega'_{j}}
\end{equation*}
avec $\Omega_{i},\ \Omega'_{j} \in \seriesfact$ telles que
$\Omega_{i}(+\infty)$ et $\Omega'_{j}(+\infty)$ soient non nulles
et $\alpha_{i}, \alpha'_{j} \in \rinteger.$ On en déduit qu'il existe $\alpha,\ \beta \in \rinteger$ tels que
$x^{\alpha}log^{\beta}(x)y$ tende vers une limite non nulle $a \in
\complex^{\ast}$ pour $x\longrightarrow +\infty.$ D'autre part, l'équation
(\ref{equation du theo de forme normale}) admet une solution de la forme $Fe_{d}^+,$
avec $F \in \seriesfact$ tangente à $I$ en $+\infty$ (on prend la
solution canonique du système fuchsien non résonnant). Posons :
\begin{equation*}
  \omega= (F e_{d}^+)^{-1} y \in Gl_n(\mathcal{M}(\complex / \rinteger)).
\end{equation*}
Ainsi $x^{\alpha}log^{\beta}(x)\omega F
  e_{d}^+$ tend vers une limite non nulle $a \in
\complex^{\ast}$ pour $x\longrightarrow +\infty.$ En notant que
$\frac{e_{d}^+}{x^{d}}$ vérifie la même propriété, on en déduit que
$\omega x^{d}x^{\alpha}log^{\beta}(x)$ admet une limite finie non
nulle pour $x\longrightarrow +\infty.$ Lorsque $d \notin -\alpha +
i \real,$ il est clair que c'est impossible. Dans le cas restant,
il existe $\lambda \in \real^{*}$ tel que $x=-\alpha + i \lambda$
et $\omega(x) x^{i \lambda}$ admet une limite finie non nulle pour
$x\longrightarrow +\infty.$ Cela implique que, étant donné $x_{0}$
un point réel qui n'est ni un pôle, ni un zéro de $\omega,$ la
suite $(x_{0}+n)^{i \lambda}$ admet une limite pour $n\in
\integer$ qui tend vers l'infini. C'est évidemment faux. Enfin, si $d=0,$ il
est clair que $\omega$ est constante; $Fe_{d}^+$ est alors
égal à $1.$
\end{proof}

\end{document}